\documentclass[journal]{new-aiaa}
\usepackage[utf8]{inputenc}

\usepackage{bm}
\usepackage{graphicx}
\usepackage{amsmath}
\usepackage[version=4]{mhchem}
\usepackage{siunitx}
\usepackage{longtable,tabularx}
\usepackage{gensymb}
\usepackage{mathtools}
\usepackage{float}
\usepackage{siunitx}
\usepackage{algorithm}
\usepackage{algpseudocode}
\usepackage{subcaption}
\usepackage{enumitem}
\usepackage{multirow}
\usepackage{lscape}
\usepackage{adjustbox}
\usepackage[flushleft]{threeparttable}
\usepackage{color,soul}
\usepackage{todonotes}

\title{Satellite Constellation Pattern Optimization for Complex Regional Coverage}

\author{Hang Woon Lee\footnote{Ph.D. Student and National Science Foundation Graduate Research Fellow, Daniel Guggenheim School of Aerospace Engineering.}}
\affil{Georgia Institute of Technology, Atlanta, GA, 30332}
\author{Seiichi Shimizu\footnote{Researcher, Mechatronics Department, Advanced Technology R\&D Center.}, Shoji Yoshikawa\footnote{Chief Researcher, Mechatronics Department, Advanced Technology R\&D Center.}}
\affil{Mitsubishi Electric Corporation, Amagasaki, Japan}
\author{Koki Ho\footnote{Assistant Professor, Daniel Guggenheim School of Aerospace Engineering. Member AIAA.}}
\affil{Georgia Institute of Technology, Atlanta, GA, 30332}

\begin{document}

\maketitle

\begin{abstract}
    The use of regional coverage satellite constellations is on the rise, urging the need for an optimal constellation design method for complex regional coverage. Traditional constellations are often designed for continuous global coverage, and the few existing regional constellation design methods lead to suboptimal solutions for periodically time-varying or spatially-varying regional coverage requirements. This paper introduces a new general approach to design an optimal constellation pattern that satisfies such complex regional coverage requirements. To this end, the circular convolution nature of the repeating ground track orbit and common ground track constellation is formalized. This formulation enables a scalable constellation pattern analysis for multiple target areas and with multiple sub-constellations. The formalized circular convolution relationship is first used to derive a baseline constellation pattern design method with the conventional assumption of symmetry. Next, a novel method based on binary integer linear programming is developed, which aims to optimally design a constellation pattern with the minimum number of satellites. This binary integer linear programming method is shown to achieve optimal constellation patterns for general problem settings that the baseline method cannot achieve. Five illustrative examples are analyzed to demonstrate the value of the proposed new approach.
\end{abstract}

\section*{Nomenclature}
{\renewcommand\arraystretch{1.0}
\noindent\begin{longtable*}{@{}l @{\quad=\quad} l@{}}
	$a$ & semi-major axis \\
	$\bm{b}$ & coverage timeline \\
	$c$ & coverage satisfactoriness indicator \\
	$e$ & eccentricity \\
	$\bm{f}$ & coverage requirement vector \\
	$i$ & inclination \\
	$\mathcal{J}$ & set of target point(s) \\
	$L$ & length (number of time steps) of vectors \\
	$M$ & mean anomaly \\
	$n$ & discrete-time instant \\
	$N$ & total number of satellites \\
	$N_{\text{D}}$ & number of Greenwich nodal periods \\
	$N_{\text{P}}$ & number of orbit nodal periods \\
	$p$ & semi-latus rectum \\
	$\bm{P}_\pi$ & permutation matrix \\
	$\bm{r}_{\text{g}}$ & target point position vector \\
	$\bm{r}_{\text{s}}$ & satellite position vector \\
	$R_{\oplus}$ & mean radius of the Earth \\
	$t$ & continuous-time instant \\
	$t_{\text{step}}$ & time step \\
	$T_{\text{G}}$ & nodal period of Greenwich\\
	$T_{\text{r}}$ & period of repetition \\
	$T_{\text{S}}$ & satellite nodal period \\
	$T_{\text{sim}}$ & simulation time horizon \\
	$\bm{v}$ & access profile \\
	$\bm{V}$ & access profile circulant matrix \\
	$\bm{x}$ & constellation pattern vector \\
	$\mathcal{Z}$ & set of sub-constellation(s) \\
	$\mathbb{Z}_2$ & binary integer number set \\
	$\mathbb{Z}_{>0}$ & positive integer number set \\
	$\mathbb{Z}_{\ge0}$ & non-negative integer number set \\
	$\varepsilon$ & elevation angle \\
	$\eta$ & satellite spacing constant \\
	$\lambda$ & longitude \\
	$\mu_{\oplus}$ & standard gravitational parameter of the Earth \\
	$\bm{\rho}$ & relative position vector from target point to satellite \\
	$\tau$ & period ratio \\
	$\phi$ & latitude \\
	$\omega$ & argument of perigee \\
	$\omega_{\oplus}$ & rotation rate of the Earth \\
	$\Omega$ & right ascension of the ascending node \\
	$\textbf{\oe}$ & orbital elements vector \\
\multicolumn{2}{@{}l}{Subscripts}\\
$j$ & target point index \\
$k$ & satellite index \\
\multicolumn{2}{@{}l}{Superscripts}\\
$z$ & sub-constellation index \\
\multicolumn{2}{@{}l}{Abbreviations}\\
APC & access profile, constellation pattern, and coverage timeline \\
BILP & binary integer linear programming \\
ECEF & Earth-centered Earth-fixed \\
ECI & Earth-centered inertial \\
GSO/GEO & geosynchronous/stationary equatorial orbit\\
NGSO & non-geostationary orbit \\
RAAN & right ascension of the ascending node \\
RGT & repeating ground track \\
\end{longtable*}}

\section{Introduction}
Satellite constellations for regional coverage are increasingly being considered as competent business solutions in a market dominated by global-based constellation systems. Regional constellations, whose form varies from being standalone to augmenting existing space-borne systems, provide flexible solutions to stakeholders as a means of circumventing geopolitical, economic, and/or technical issues associated with global constellation systems. Examples of such regional constellation systems are the Indian Regional Navigation Satellite System (IRNSS) \cite{irnss} and the Quasi-Zenith Satellite System \cite{qzss}.

Unlike global coverage constellations, regional coverage constellations solely focus on the coverage over a local region and therefore generally require a smaller number of satellites in the system to achieve the same performance per area metric compared to global-coverage constellations. This leads to a significantly reduced system cost as the total life-cycle cost of the system depends on the number of satellites \cite{diekelman1998}. The reduced system cost allows for a tolerable risk of failure and facilitates a shorter payback period. These properties allow regional constellation systems to swiftly react to uncertainties arising from market demand and/or administrative issues. Research has also shown that a flexible option to treat a regional constellation system as part of a larger staged deployment process can be beneficial when market uncertainties are present \cite{lee2018}.

Various space systems have been designed for regional coverage. Although the most classical regional coverage method is to use geo-synchronous/-stationary equatorial (GSO/GEO) orbits, non-geostationary orbit (NGSO) systems are deemed to provide better performance for many mission-critical attributes such as latency and launch cost. Traditional constellation design methods have investigated the problems with relatively simple coverage criteria, such as satisfying an $f$-fold continuous coverage requirement (e.g., single-fold, double-fold, etc.) or minimizing the maximum revisit time gap over an area. However, the problems with complex coverage requirements that are periodically time-varying and spatially-varying have not been explored. Examples of such coverage requirements are (1) the increased communication service needs during the daytime; and (2) the increased service needs over urban/sub-urban areas for reliable access \cite{lutz2012}. The design process to generate the optimal constellation for such complex coverage requirements involves determining (1) the common orbital characteristics and (2) a constellation pattern. While conventional constellation design methods often assume a symmetric pattern (e.g., Walker constellations \cite{walker1970,walker1977,walker1984}) and optimize the common orbital characteristics (e.g., altitude, inclination) with that assumption, the large design space of asymmetric constellation patterns is often missed despite its importance particularly for complex time-varying and spatially-varying coverage requirements. Furthermore, it is reasonable to assume that such a regional coverage constellation system can constitute multiple sub-constellations, each with different orbital characteristics, as demonstrated in the case of IRNSS \cite{irnss}; however, the concurrent design of multiple sub-constellation patterns using NGSOs requires a sophisticated optimization approach. Such a topic has been scarcely studied and remains an open question. Given this background, a research question of interest arises: ``\textit{How do we design a constellation pattern (for multiple sub-constellations if needed) that is optimized (i.e., with the minimum number of satellites) for a periodically time-varying and spatially-varying demand over the regional area(s) of interest?}'' This paper seeks to address this question by constructing an optimal constellation pattern design approach for complex regional coverage. The resulting rigorous constellation pattern design approach can be integrated with existing orbital characteristics design methods and launch/mission constraints to optimize future satellite constellation design.

The contribution of this paper is as follows. First, the discovery of a \textit{circular convolution phenomenon} between a seed satellite access profile, a constellation pattern vector, and a coverage timeline is formalized in this research. The resulting formulation is referred to as the APC decomposition, following the acronyms of the seed satellite Access profile, constellation Pattern, and Coverage timeline; each of these concepts is introduced in detail later in this paper. We derive a linear formulation that enables us to design a constellation pattern for a system of multiple sub-constellations for multiple regions. This formulation provides a foundation for general methods introduced herein. 
Second, we extend the traditional definition of a time-independent $f$-fold coverage requirement (e.g., single-fold, double-fold, etc.) to a time-dependent $f[n]$-fold coverage requirement, where $n$ is a discrete-time instant, such that periodically time-varying coverage demands can be handled optimally in the constellation design. By applying this idea to multiple target points, this approach is further extended to the case with time-varying and spatially-varying coverage requirements.
Finally, we develop a general method based on binary integer linear programming (BILP) that finds the optimal satellite constellation pattern for complex regional coverage, and, if needed, for multiple sub-constellations concurrently. This core concept enables users to explore the hidden design space by breaking the symmetry in the constellation design. The developed constellation pattern design approach is demonstrated with a series of case study examples.

The rest of the paper is organized as follows. Section~\ref{sec:lit_review} provides a summary of the key literature relevant to this research. Section~\ref{sec:sat_con_gen} provides an overview of the constellation model used in this paper. Section~\ref{sec:circularconvolutionformulation} introduces the ideas behind the developed approach, including the circular convolution formulation of the problem and its pertinent definitions. Section~\ref{sec:methodology} then introduces two methods based on this formulation: the baseline quasi-symmetric and the novel BILP methods. The developed methods are applied to various illustrative examples in Section~\ref{sec:illustrative_examples} for demonstration. Section~\ref{sec:conclusion} then concludes this paper.

\section{Literature Review} \label{sec:lit_review}
This section reviews the major literature relevant to this study. Traditional satellite constellation design methods have focused on minimizing the number of satellites while providing continuous coverage over a large area of interest such as the globe or latitudinally-bounded zones. Classical methods such as the streets of coverage \cite{luders1961, luders1974, beste1978, rider1986}, Walker and Rosette constellations \cite{walker1970,walker1977,walker1984,ballard1980}, and the tetrahedron elliptical constellation \cite{draim1987} leveraged a geometric approach to exhibit a symmetry in the constellation pattern, where satellites are uniformly and symmetrically arranged based on a predetermined phasing rule. The symmetry in the constellation pattern provides a foundation for a complete design space analysis due to finite variability \cite{wertz2001} or for an analytical solution. However, this usually leads to redundant coverage overlaps and therefore may not produce an optimal constellation design in terms of the number of satellites over a bounded local region.

There are several prior studies that specifically dealt with the design of regional coverage constellations. By fully utilizing the characteristics of the repeating ground track orbits, Hanson et al. \cite{hanson1992designing} and Ma and Hsu \cite{ma1997} utilized the timeline meshing method to generate the optimal constellation with respect to minimizing the maximum time gap at the minimum possible inclination. Similarly, Pontani and Teofilatto extended the characteristics of the repeating ground track by searching for allowable time delays with respect to minimizing the gap or maximizing coverage \cite{pontani2007}. In addition, Crossley and Williams used metaheuristics methods to design a satellite constellation to minimize the maximum revisit time \cite{crossley2000}. Although these regional constellation design algorithms show promising ability to produce asymmetric configuration with respect to a single target point or a connected area, these methods are not applicable to designing a constellation system for periodically time-varying demands over multiple disjoint target points (referred to as complex coverage requirement in this paper) with multiple sub-constellations. Ulybyshev investigated a new geometric approach to generate satellite constellation designs for complex coverage \cite{ulybyshev2008}. The method demonstrates the use of the two-dimensional space and combined maps for the satellite constellation and coverage functions. Nevertheless, this method cannot be applied to asymmetric constellations. Other literature can be found in the comprehensive literature review by Dutruel-Lecohier and Mora as well as Wertz \cite{dutruel-lecohier1998,wertz2001}. Recently, Ulybyshev presented a short historical survey of satellite constellation design for continuous coverage \cite{ulybyshev2009}. However, there is no methodology that directly answers our question raised in the introduction that considers all three aspects of the regional coverage problem: (1) multiple target points, (2) complex coverage requirements, and (3) multiple sub-constellations.

In response to this background, this paper attempts to construct methods to design a satellite constellation pattern for periodically time-varying and spatially-varying coverage requirements over multiple target points, and if demanded, for multiple sub-constellations. Building upon the idea of repeating ground track orbits and common ground track constellations (e.g., Flower Constellation set theory \cite{mortari2004,mortari2008,wilkins2008}) and generalizing our prior work \cite{lee2018AAS}, we formalize the circular convolution nature of the constellation pattern design problem and derive two methods for it: (1) the baseline and rather traditional quasi-symmetric method; and (2) the more general and novel BILP method. The developed approach can design the constellation pattern that satisfies the complex coverage requirements of multiple target points with the minimum number of satellites possible exploring both symmetric and asymmetric patterns.

\section{Satellite Constellation Model}
\label{sec:sat_con_gen}
This section introduces the ideas and assumptions on the satellite constellation model that the proposed approach builds upon, including the repeating ground track orbit and the common ground track constellation. 
\subsection{Repeating Ground Track Orbit}
A ground track of a satellite is defined as a trace of its sub-satellite points on the surface of the Earth. In this paper, we utilize a repeating ground track (RGT) orbit as a basis for the orbital design of the constellation, which allows a ground track of a satellite to repeat exactly and periodically. This type of orbit has been shown to provide better coverage performance than the non-repeating ground track orbits with fewer satellites for regional coverage \cite{hanson1992designing}. Considering the Earth-centered Earth-fixed (ECEF) frame, an RGT orbit is achieved when the nodal period of the orbit $T_{\text{S}}$ (the time interval between two consecutive crossings of the orbit ascending node by a satellite) is a rational multiple of the nodal period of Greenwich $T_{\text{G}}$ (the time interval between two consecutive crossings of the orbit ascending node line by the prime meridian):

\begin{equation}
	\label{eq:period_of_repetition}
	T_{\text{r}} = N_{\text{P}}T_{\text{S}}=N_{\text{D}}T_{\text{G}}
\end{equation}

{\parindent0pt
	where $T_{\text{r}}$ represents the period of repetition. Eq.~\eqref{eq:period_of_repetition} implies that a satellite on an RGT orbit makes $N_{\text{P}}$ number of revolutions in $N_{\text{D}}$ number of nodal periods of Greenwich \cite{mortari2004,vtipil2012determining}. $N_{\text{P}}$ and $N_{\text{D}}$ are positive integer numbers.
}

Considering the $J_2$ perturbation effect, the nodal period of the satellite orbit $T_{\text{S}}$ and the nodal period of Greenwich $T_{\text{G}}$ are given in Eqs.~\eqref{eq:t_s_t_g}:

\begin{subequations}
\label{eq:t_s_t_g}
\begin{align}
T_{\text{S}} = \frac{2\pi}{\dot{\omega}+\dot{M}} \\
T_{\text{G}} = \frac{2\pi}{\omega_{\oplus}-\dot{\Omega}} \label{eq:tg}
\end{align}
\end{subequations}

{\parindent0pt
	where $\omega_{\oplus}$ is the rotation rate of the Earth, $\dot{\omega}$ is the rate of change in the argument of perigee due to perturbations, $\dot{\Omega}$ is the rate of nodal regression of a satellite's orbit, and $\dot{M}$ is the rate of change in the mean anomaly due to nominal motion and perturbations. The perturbed orbital elements in Eqs.~\eqref{eq:t_s_t_g} are:
}

\begin{eqnarray}
\dot{\omega} = \frac{3}{2}J_2\Bigg(\frac{R_{\oplus}}{p}\Bigg)^2 \sqrt{\frac{\mu_{\oplus}}{a^3}}\Bigg[2-\frac{5}{2}\sin^2i\Bigg] \\
\dot{\Omega} = -\frac{3}{2}J_2\Bigg(\frac{R_{\oplus}}{p}\Bigg)^2\sqrt{\frac{\mu_{\oplus}}{a^3}}\cos i \\
\dot{M} = \sqrt{\frac{\mu_{\oplus}}{a^3}}\Bigg[1-\frac{3}{2}J_2\Bigg(\frac{R_{\oplus}}{p}\Bigg)^2\sqrt{1-e^2}\Bigg(\frac{3}{2}\sin^2i-1\Bigg)\Bigg]
\end{eqnarray}

{\parindent0pt
	where $R_{\oplus}=\SI{6378.14}{km}$ is the mean radius of the Earth, $p=a(1-e^2)$ is the semi-latus rectum, $\mu_{\oplus}=\SI{398600.44}{km^3 s^{-2}}$ is the standard gravitational parameter of the Earth, and $J_2=0.00108263$  is the zonal harmonic coefficient due to the equatorial bulge of the Earth \cite{wertz2001}.
}

A period ratio $\tau$ is defined as a ratio of $N_{\text{P}}/N_{\text{D}}$ and further can be deduced based on the perturbed orbital elements:
\begin{equation}
	\label{eq:tau}
	\tau=\frac{N_{\text{P}}}{N_{\text{D}}}=\frac{T_{\text{G}}}{T_{\text{S}}}=\frac{\dot{\omega}+\dot{M}}{\omega_{\oplus}-\dot{\Omega}}
\end{equation}
The period ratio is used to identify a unique RGT orbit out of an $N_{\text{P}}$ and $N_{\text{D}}$ pair \cite{mortari2004}. That is, a satellite orbit with $\tau = 10/2$ and a satellite orbit with $\tau=5/1$ both of which share an identical orbit and a ground track.

The semi-major axis $a$ of an RGT orbit can be derived using the Newton-Raphson method presented by Bruccoleri for a given set of $N_{\text{P}}$, $N_{\text{D}}$, $e$, and $i$ \cite{bruccoleri2007}. Because the semi-major axis is a function of $\tau$, $e$, and $i$ (i.e., $a=a(\tau,e,i)$), we shall utilize the period ratio $\tau=N_{\text{P}}/N_{\text{D}}$ as an independent orbital variable instead of the semi-major axis $a$. Henceforth, this paper utilizes an RGT orbital elements vector, $\textbf{\oe}=[\tau,e,i,\omega,\Omega,M]^T$, to fully define an RGT orbit of a satellite. We assume the utilization of satellite maneuvers to correct and maintain an identical ground track throughout the satellite lifetime, negating perturbation effects other than the $J_2$ effect. Note that the right ascension of the ascending node (RAAN) $\Omega$ and the mean anomaly $M$ in the RGT orbital elements vector indicate the initial values in reference to a given epoch.

\subsection{Common Ground Track Constellation}
\begin{figure}[H]
	\centering\includegraphics[width=0.55\textwidth]{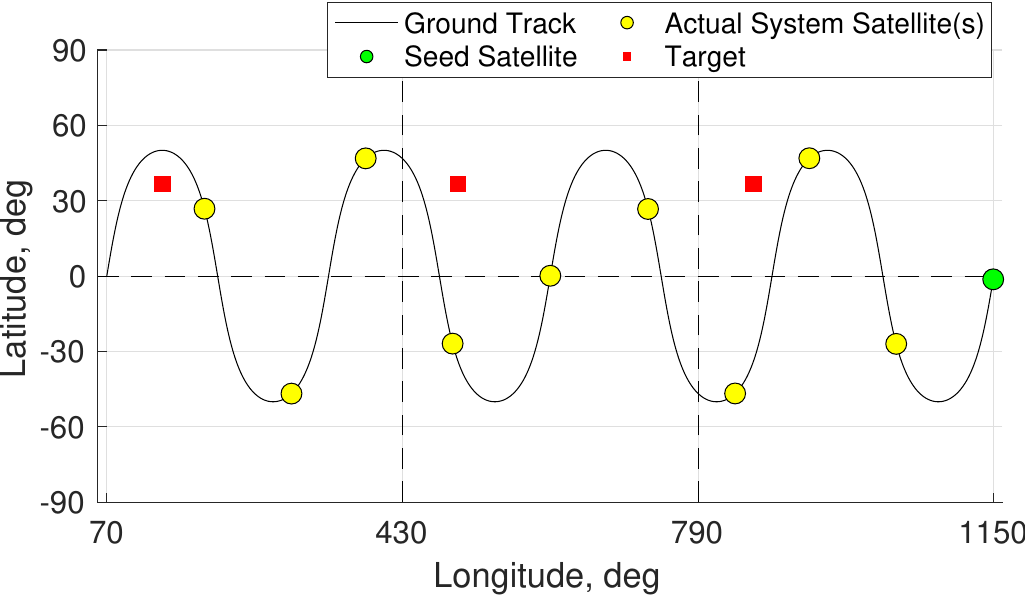}
	\caption{Illustration of a common ground track constellation in an expanded ground track view}
	\label{fig:cgc}
\end{figure}

This paper considers a constellation pattern where all satellites in the constellation are systematically generated such that their ground tracks overlap to create a single common ground track. In this paper, we refer to this type of constellation as a \textit{common ground track constellation}. (If there are multiple sub-constellations, each sub-constellation has its own common ground track.) Fig.~\ref{fig:cgc} illustrates an example of arbitrarily defined 9-satellite common ground track constellation; its system satellites, depicted in yellow circles, are placed along a common ground track. The definitions of the terms used in Fig.~\ref{fig:cgc} are discussed in Section~\ref{sec:circularconvolutionformulation}. For more information about the expanded ground track view, refer to Appendix~A.

A common ground track constellation has relationships with several constellation design theories. For example, when certain conditions are satisfied (e.g., symmetric distribution and $N_\text{D}=1$), the common ground track constellations utilizing circular RGT orbits can be expressed as $i:N/N/(N-N_{\text{P}})$ when $N_{\text{D}}=1$ \cite{wu2006design}, following the standard Walker notation $i: N/P/F$. Here, $N$ is the total number of satellites in the system, $P$ is the number of orbital planes, and $F$ is the Walker phasing factor.

A common ground track constellation with RGT orbits is a common assumption used in the literature such as the original Flower Constellation theory \cite{mortari2004}. The Flower Constellation is defined as a set of $N$ satellites following the same (closed) trajectory with respect to a rotating frame. For this paper, the ECEF frame is considered. Ref.~\cite{avendano2013} introduces the three conditions to construct a Flower Constellation as follows:

\begin{enumerate}
	\item The orbital period of each satellite is a rational multiple of the period of the rotating frame. 
	\item The semi-major axis $a$, eccentricity $e$, inclination $i$, and argument of perigee $\omega$ are identical for all the satellite orbits. 
	\item The right ascension of the ascending node $\Omega_k$ and the mean anomaly $M_k$ of each satellite ($k=1,...,N$) satisfy:
	\begin{equation}
	\label{eq:flower}
	N_{\text{P}} \Omega_k+N_{\text{D}} M_k = \text{constant} \ \text{mod} \ (2\pi)
	\end{equation}
\end{enumerate}
This paper utilizes the above three conditions of the original Flower constellation set theory as a basis for constellation generations. Furthermore, we restrict satellite orbits to be either circular or critically-inclined elliptic ($i=63.4\degree$ or $116.6\degree$). This is because, in engineering practice, non-critically-inclined elliptic orbits are generally avoided for periodic coverage requirements due to heavy orbital maintenance costs incurred by negating the precession of the argument of perigee.

\section{Circular Convolution Formulation} \label{sec:circularconvolutionformulation}
Building upon the satellite constellation model in the previous section, this section introduces the main ideas behind the methods developed in this paper, including the definitions and concepts of the access profile, coverage, and constellation pattern representation, as well as the mathematical representation of the circular convolution phenomenon. 

The derivation of the circular convolution phenomenon utilizes time discretization. One underlying assumption is that, to satisfy the periodically time-varying coverage requirements, the repeat period of the RGT orbit can be chosen such that it is a rational multiple of the repeat period of the coverage requirement. This implies that we can discretize both of these repeat periods by a common time step length $t_\text{step}$. The least common multiple of the numbers of time steps for these two repeat periods would be the number of time steps for the simulation time horizon length $L$ needed to evaluate the coverage.
If there are multiple target points with different repeat periods for their coverage requirements, assuming that their repeat periods can each be represented as an integer number of time steps with the common interval $t_\text{step}$, then we can use the least common multiple of these time steps as the repeat period of the ``overall'' coverage requirement. The circular convolution formulation and its associated properties are defined over this discretized $L$-step simulation time horizon length.

For simplicity, in this paper, we consider the case in which the repeat period of RGT is an integer multiple of the repeat period of the coverage requirement; in this case, $T_\text{sim}=T_\text{r}=Lt_\text{step}$, where $T_\text{sim}$ is the length of the simulation time horizon and $T_\text{r}$ is the repeat period of the orbit. This case can be easily generalized to the above more general case. Note that the uniformly continuous coverage case can be treated as a special case, where the repeat period of the orbit (and thus the simulation time horizon) can be arbitrarily chosen.

\subsection{Access Profile} \label{sec:accessprofile}
The relative position vector $\bm{\rho}$ pointing from a ground target point to a satellite is defined as:

\begin{equation}
\bm{\rho}=\bm{r}_{\text{s}}-\bm{r}_{\text{g}}
\end{equation}

\begin{figure}[H]
	\centering\includegraphics[width=0.23\textwidth]{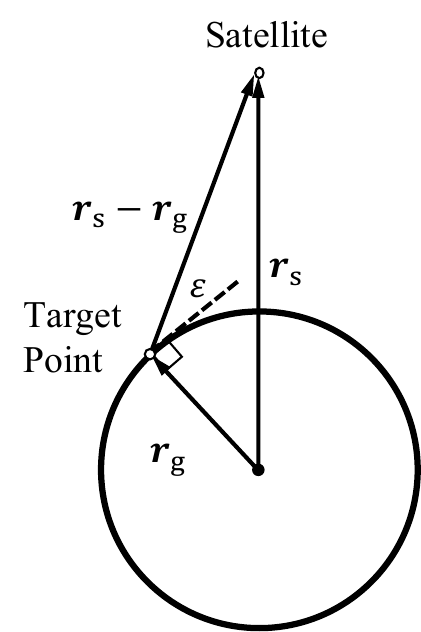}
	\caption{Satellite, target point, and elevation angle relationship}
	\label{fig:vector}
\end{figure}

{\parindent0pt
	where $\bm{r}_{\text{s}}$ is a satellite position vector from the center of the Earth and $\bm{r}_{\text{g}}$ is a target point position vector from the center of the Earth. Fig.~\ref{fig:vector} illustrates this relationship.
}

An elevation angle $\varepsilon$ of a satellite seen from a ground target point is defined as:

\begin{equation}
\varepsilon = \sin^{-1} \Bigg(\frac{\bm{r}_{\text{g}}\cdot\bm{\rho}}{\lVert\bm{r}_{\text{g}}\rVert \lVert\bm{\rho}\rVert}\Bigg) = \sin^{-1} (\hat{\bm{r}}_{\text{g}} \cdot \hat{\bm{\rho}})
\end{equation}

{\parindent0pt
	where $\lVert \cdot \rVert$ is the Euclidean norm.
}

\begin{figure}[H]
	\centering\includegraphics[width=0.6\textwidth]{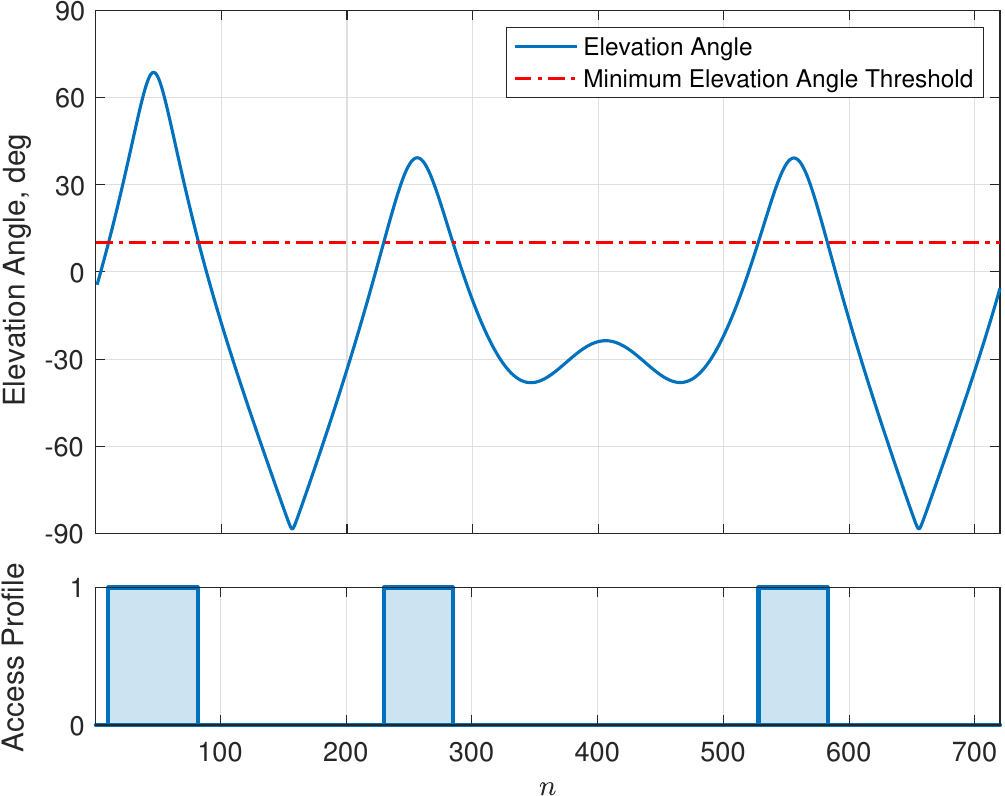}
	\caption{Sample illustration of a satellite's elevation angle viewed from a ground point and corresponding access profile}
	\label{fig:elevation}
\end{figure}

Because the dot product between the unit target point position vector $\hat{\bm{r}}_{\text{g}}$ and the unit relative position vector $\hat{\bm{\rho}}$ continues to change due to the rotation of the Earth and the motion of a satellite, the elevation angle is, therefore, a function of time, $\varepsilon=\varepsilon(t)$. An example of a typical NGSO satellite elevation angle function is shown in the upper part of Fig.~\ref{fig:elevation}. When the elevation angle of a satellite is above the minimum elevation angle threshold $\varepsilon_{\text{min}}$, which is determined by the mission requirement \cite{lutz2012}, the satellite is said to be visible from or \textit{to have access to the target point}. Since the periods when the satellite has access to the ground target point are of particular interest, we convert the elevation angle function into an \textit{access profile} (or a visibility profile in some literature), which is a binary vector that indicates either access, 1, or no access, 0, at each time instant. The access profile is visualized in the lower part of Fig.~\ref{fig:elevation}. This paper utilizes a sampling method to generate an access profile. Note that access profiles can be derived in different ways \cite{chylla1992,alfano1992,han2017}.

The continuous-time elevation angle function $\varepsilon(t)$ is sampled at every time step of $t_{\text{step}}$ to create a discrete-time elevation angle function $\varepsilon[n]$ with length $L$. As mentioned earlier, $L$ is the number of time steps of the simulation horizon, i.e., $T_{\text{sim}}=Lt_{\text{step}}$, where $T_{\text{sim}}$ is the simulation time horizon (which is assumed to be equal to the RGT repeat period $T_\text{r}$ in this paper for simplicity as discussed above). The access profile $\bm{v}_{k,j}\in\mathbb{Z}^{L}_{2}$ between the $k$th satellite and the $j$th target point stores boolean information about the satellite access (or visibility) state at each discrete-time instant $n\in\{0,...,L-1\}$. Therefore, each element of the access profile is:

\begin{equation}
\label{eq:access}
v_{k,j}[n]\triangleq\begin{cases}
1, &\text{if $\varepsilon_{k,j}[n]\geq\varepsilon_{k,j,\text{min}}[n]$} \\
0, &\text{otherwise}
\end{cases}
\end{equation}

{\parindent0pt
	where $n$ is the discrete-time instant and $\mathcal{J}$ is the set of target points. Throughout this paper, vectors are represented in italic boldface (e.g., $\bm{v}_{k,j}$) and their elements are represented in brackets (e.g., $v_{k,j}[n]$). To make the notation consistent with the circular convolution method from the digital signal processing community, the vector index representing the discrete-time instant $n$ is set to take the range of $[0,L-1]$.
}

It is important to note a condition in Eq.~\eqref{eq:access}: there must exist at least one access interval for a given satellite-target for the methods introduced in this paper to function; simply stated, the access profile shall be a non-zero vector. The methods introduced in the following sections are constructed based on the assumption that the access profile is a non-zero vector.

One can interpret the generalized minimum elevation angle $\varepsilon_{k,j,\text{min}}[n]$ in Eq.~\eqref{eq:access} as the minimum elevation angle threshold imposed on an access between a satellite $k$ and a target point $j$ at discrete-time instant $n$. This paper assumes that all satellites have a common generalized minimum elevation angle:

\begin{equation}
\label{eq:epsilon_min}
\varepsilon_{k,j,\text{min}}[n]=\varepsilon_{j,\text{min}}[n], \ \ \ \forall k \in\{1,...,N\}
\end{equation}

When designing a satellite constellation for regional coverage, a constellation must be spatially and temporally referenced relative to the target point and the epoch. A hypothetical satellite that conveys referenced orbital information, $\textbf{\oe}_0=[\tau,e,i,\omega,\Omega_0,M_0]^T$, for the constellation is defined as the \textit{seed satellite}\footnote{The term seed satellite is credited to the software Systems Tool Kit (STK) \cite{stk}.} and the corresponding $\textbf{\oe}_0$ as the \textit{seed satellite orbital elements vector}. The actual satellites inherit the common orbital characteristics defined in this seed satellite elements vector, but they independently hold ($\Omega_k,M_k$) pairs that are determined by Eq.~\eqref{eq:flower}, resulting in the orbital elements vector for each satellite of $\textbf{\oe}_k=[\tau,e,i,\omega,\Omega_k,M_k]^T$ where $k$ is an index of a satellite ($\Omega$ and $M$ are initial values referenced to a given epoch; the subscripts refer to the index of a corresponding satellite). Note that it is not required to have an actual satellite at the seed satellite position; the seed satellite orbital elements are used as a reference to define the actual satellites in the system.

Let us recall the main assumptions considered thus far: (1) all satellites are placed on a common repeating ground track constellation as shown in Fig.~\ref{fig:cgc}; and (2) all access between a target point $j$ and every member satellite in a given constellation are constrained to the same minimum elevation angle threshold as shown in Eq.~\eqref{eq:epsilon_min}. Such assumptions enable us to utilize a powerful property, a \textit{cyclic property}, in which all access profiles of the member satellites in a given constellation are \textit{identical, but circularly shifted}. Therefore, any access profile $\bm{v}_{k,j}$ between the $k$th satellite and the $j$th target point can be represented as a \textit{circularly shifted seed satellite access profile} $\bm{v}_{0,j}$:

\begin{equation}
\label{eq:bb}
v_{k,j}[n]=\bm{P}_\pi^{n_k}v_{0,j}[n]
\end{equation}

{\parindent0pt
	where $\bm{P}_\pi$ is a permutation matrix with the dimension $(L \times L)$ as shown in Eq.~\eqref{eq:permutation} and $n_k$ is the index representing its (temporal) location of the $k$th satellite with respect to the seed satellite along the common ground track. 
}
	
\begin{equation}
\label{eq:permutation}
\bm{P}_\pi=
\renewcommand{\arraystretch}{0.8}
\begin{bmatrix}
0 & 0 & 0 & \cdots &1 \\
1 & 0 & 0 & \cdots & 0 \\
0 & 1 & 0 & \ddots & \vdots \\
\vdots & \vdots & \ddots & \ddots & 0 \\
0 & 0 & \cdots & 1 & 0
\end{bmatrix}
\end{equation}
The formal definition and the physical interpretation of $n_k$ are explained in Section~\ref{sec:x}.

\subsection{Coverage Timeline and Coverage Requirement}
\begin{figure}[H]
	\centering\includegraphics[width=0.3\textwidth]{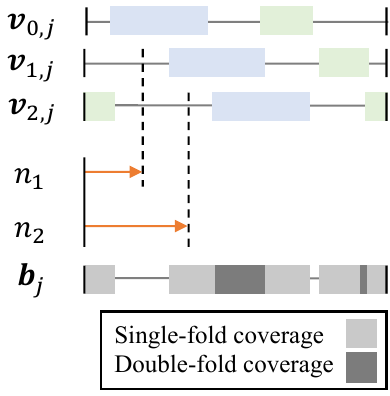}
	\caption{Illustration of shifts of access profiles (2-satellite system); notice that the seed satellite access profile ($\bm v_{0,j}$) is not part of the coverage timeline}
	\label{fig:shift}
\end{figure}

Because there are multiple satellites in the constellation system, the access profiles must be meshed together to create a \textit{coverage timeline} over a target point. Hence, a coverage timeline $\bm{b}_j\in\mathbb{Z}^{L}_{\ge0}$ is an access profile between multiple satellites and the $j$th target point; it stores information about the number of satellites in view at each discrete-time instant $n$. This is illustrated in Fig.~\ref{fig:shift}. As before, $n_1$ and $n_2$ are the indices that represent the temporal locations of the first ($k=1$) and second ($k=2$) satellite with respect to the seed satellite ($k=0$), respectively. It is important to point out that because the seed satellite is hypothetical, its access profile is not considered in the coverage timeline. Eq.~\eqref{eq:b} provides a mathematical definition of the coverage timeline:

\begin{equation}
\label{eq:b}
b_j[n]=\sum_{k=1}^{N} v_{k,j}[n]
\end{equation}
Note that the coverage timeline is not a binary vector, but instead, it is a non-negative integer vector.

Next, we define the \textit{coverage requirement}. A coverage requirement $\bm{f}_j\in\mathbb{Z}^{L}_{\ge0}$ is a vector of non-negative integers that is created by a user per mission requirement. It is important to distinguish the difference between the coverage timeline $\bm{b}_j$ and the coverage requirement $\bm{f}_j$. The coverage timeline is a coverage performance or a state of a constellation system, whereas the coverage requirement indicates what a constellation system shall achieve. For example, in order for a constellation system to achieve $f$-fold continuous coverage, the coverage timeline must be greater than or equal to the coverage requirement, that is, at least $f$ satellite(s) must have access to or be visible by the target point throughout the simulation time horizon. The coverage satisfactoriness indicator $c_j$ indicates the coverage requirement satisfactoriness of the coverage timeline over a target point $j$.

\begin{equation}
\label{eq:c_j}
c_j\triangleq\begin{cases}
1, &\text{if} \ b_j[n] \ge f_j[n], \ \ \forall n \in \{0,...,L-1\} \\
0, &\text{otherwise}
\end{cases}
\end{equation}

If an area of interest consists of multiple target points (e.g., due to area grid discretization), the coverage is satisfactory if all target points are satisfactorily covered. Extending Eq.~\eqref{eq:c_j}, the satisfactory condition of the coverage over all target points in a set $\mathcal{J}$ can be expressed as:

\begin{equation}
\label{eq:c}
c_{\mathcal{J}}\triangleq\begin{cases}
1, &\text{if} \ c_j=1, \ \ \forall j\in \mathcal{J} \\
0, &\text{otherwise}
\end{cases}
\end{equation}

{\parindent0pt
	where $\mathcal{J}$ is a set of target points. Thus, designers of the constellation system must aim to satisfy all coverage requirements on every target point as each target point may impose its own unique coverage requirement.
}

\subsection{Constellation Pattern Vector} \label{sec:x}
\begin{figure}[H]
	\centering\includegraphics[width=0.7\textwidth]{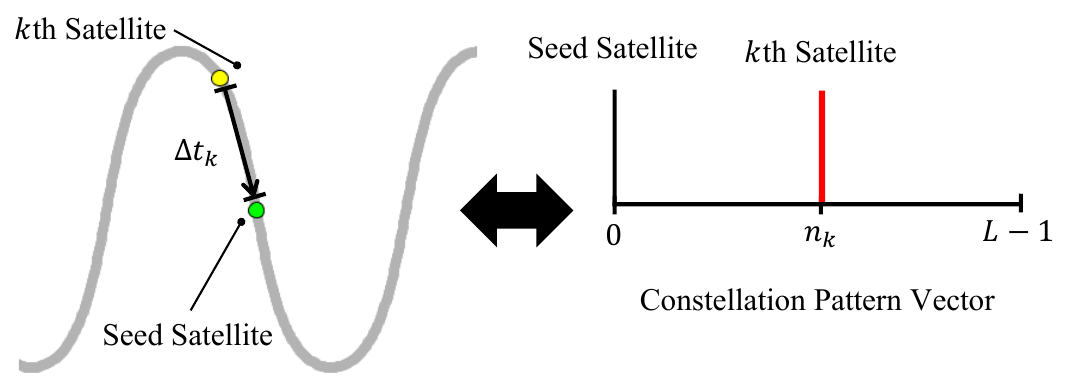}
	\caption{Illustration of a satellite time shift and its representation in the constellation pattern vector form; the direction of the motion of satellites is indicated by the arrow on the left-hand side of the figure}
	\label{fig:pattern}
\end{figure}

We express the time shifts of satellites along the ground track with respect to the seed satellite in a discrete-time binary sequence $\bm{x}\in\mathbb{Z}^{L}_{2}$ and refer to it as the \textit{constellation pattern vector}.

\begin{equation}
\label{eq:x}
x[n]\triangleq\begin{cases}
1, &\text{if } n=n_k \\
0, &\text{otherwise}
\end{cases}
\end{equation}
The temporal location index, $n_k$, can be interpreted as the time-delay index for the $k$th satellite. This is because the $k$th satellite that is delayed behind the seed satellite by the time difference of $\Delta t_k=t_{\text{step}}n_k$ over the common ground track can be equivalently shown as a unit impulse at time instant $n=n_k$ on a constellation pattern vector $\bm{x}$. This is illustrated in Fig.~\ref{fig:pattern}. The left-hand side of the figure shows a snapshot of an arbitrary constellation system: a seed satellite depicted as the green circle and an arbitrary $k$th satellite depicted as the yellow circle in an expanded ground track view. The $k$th satellite is positioned behind the seed satellite in a moving direction by the time unit of $\Delta t_k$. That is, the $k$th satellite will occupy the current position of the seed satellite $\Delta t_k$ time units later (i.e., $n_k$ time steps later). The equivalent representation in the constellation pattern vector form is shown on the right-hand side of the figure. In this case, the position of the $k$th satellite is represented as a red impulse, which represents the time delay with respect to the seed satellite.

From Eq.~\eqref{eq:x} and because $L$ is assumed to cover exactly one repeat period of the RGT, we can deduce the total number of satellites in the constellation from the constellation pattern vector as:

\begin{equation}
\label{eq:n_sat}
N=\sum_{n=0}^{L-1} x[n]
\end{equation}

\subsection{Circular Convolution Phenomenon} \label{sec:circularconvolution}
The discrete-time sequences, $\bm{v}_{0,j}$, $\bm{x}$, and $\bm{b}_j$, defined in the previous sections have a finite periodic length of $L$ due to the cyclic property of the closed relative ground track assumption. Note that, as mentioned earlier, this length of the vectors is the total number of time steps for the simulation time horizon.

A discrete circular convolution operation between the seed satellite access profile $\bm{v}_{0,j}$ and the constellation pattern vector $\bm{x}$ produces a coverage timeline $\bm{b}_j$:

\begin{equation}
\label{eq:circular}
\begin{split}
b_j[n]&=v_{0,j}[n] \circledast x[n] = \sum_{m=0}^{L-1} v_{0,j}[m] x[(n-m) \bmod L] \\
&=x[n] \circledast v_{0,j}[n] = \sum_{m=0}^{L-1} x[m] v_{0,j}[(n-m) \bmod L]
\end{split}
\end{equation}
where $\circledast$ represents a circular convolution operator. (Note that the circular convolution is commutative.) Or equivalently, this equation can be written as:

\begin{equation}
\label{eq:circulantX}
\bm{V}_{0,j}\bm{x}=\bm{b}_j
\end{equation}

{\parindent0pt
	where $\bm{V}_{0,j}\in\mathbb{Z}^{L\times L}_{2}$ is a seed satellite access profile circulant matrix that is fully specified by a seed satellite access profile $\bm{v}_{0,j}$. Note that a circulant matrix is a special form of a Toeplitz matrix \cite{gray2006toeplitz}; each entry of the matrix $[\alpha,\beta]$ is defined as:
}

\begin{equation}
\label{eq:vdef}
V_{0,j}[\alpha,\beta] = v_{0,j}[(\alpha-\beta) \bmod L]
\end{equation}

{\parindent0pt
where $\alpha$ and $\beta$ are the row and column indices, respectively, for $\alpha,\beta\in\{0,1,\cdots,L-1\}$. More information about the circular convolution is referred to Ref.~\cite{dsp}. The derivation of the circular convolution relationship (Eq.~\eqref{eq:circulantX}) from Eq.~\eqref{eq:b} is described in Appendix~B.}

\begin{figure}[H]
	\centering
	\begin{subfigure}[h]{0.495\linewidth}
	    \centering
		\includegraphics[height=6.5cm]{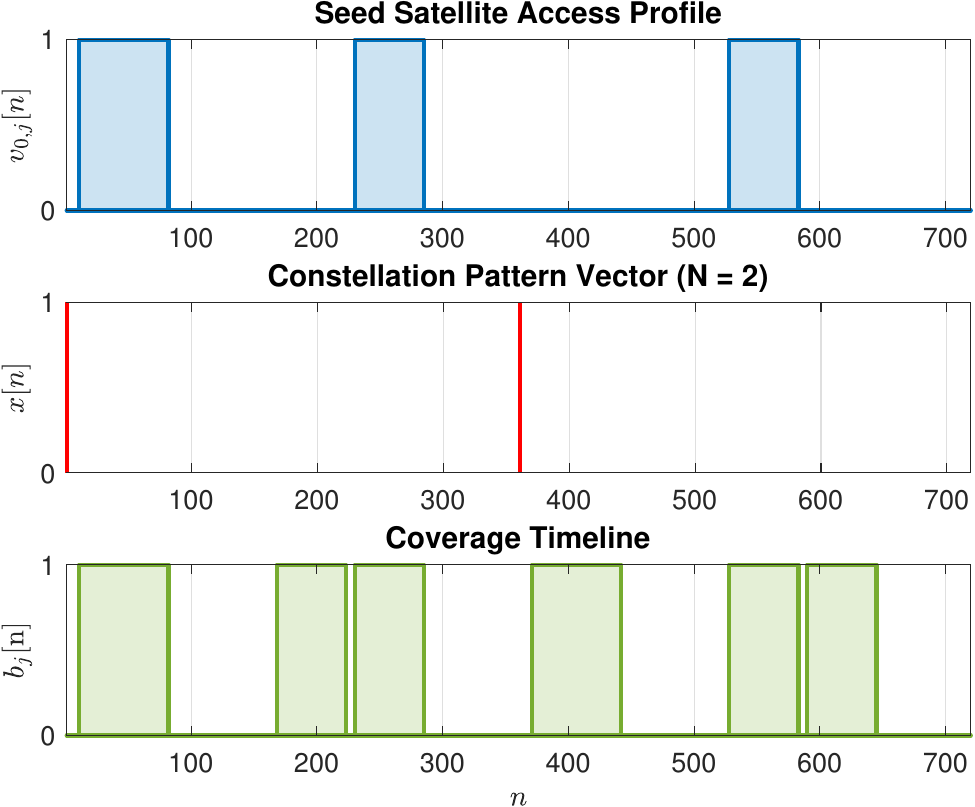}
		\caption{The APC decomposition}
		\label{fig:example:a}
	\end{subfigure}
	\begin{subfigure}[h]{0.495\linewidth}
	    \centering
		\includegraphics[height=6.5cm]{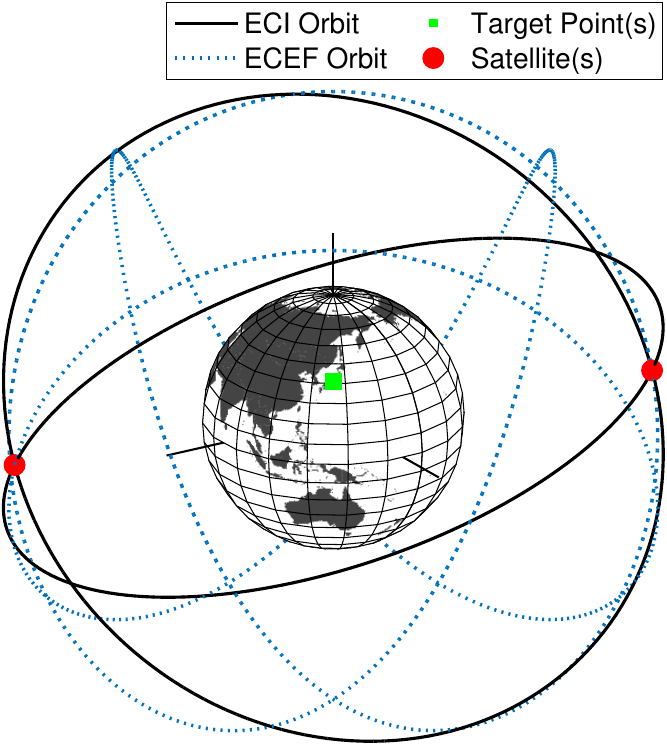}
		\caption{The corresponding configuration in the ECI and ECEF frame at $n=0$}
		\label{fig:example:b}
	\end{subfigure}
	\caption{APC decomposition and its equivalent constellation representation in 3D space}
	\label{fig:example}
\end{figure}

To illustrate this relationship, consider a system with $\textbf{\oe}_0=[4/1,0,50\degree,0\degree,350.2\degree,0\degree]^T$ (J2000) and uniformly spaced $N=2$ satellites. The corresponding seed satellite access profile observed from a target $\mathcal{J}=\{(\phi=36.7\degree \text{N},\lambda=137.48\degree \text{E})\}$ (a point) with $\varepsilon_{\text{min}}=10\degree$ is shown in the top part of Fig.~\ref{fig:example:a}. In this example, the length of vectors is set to $L=720$ such that the corresponding time step is approximately \SI{120}{s}. The constellation pattern vector, shown in the middle part of Fig.~\ref{fig:example:a}, has two unit impulses at $n=0$ and $n=360$ to represent the temporal locations of two satellites with respect to the seed satellite. The equivalent orbital elements vectors for these satellites are (refer to Section~\ref{sec:post} for the derivation):
\begin{equation*}
\label{eq:equivalent_vectors}
\begin{split}
\textbf{\oe}_1&=[4/1,0,50\degree,0\degree,350.2\degree,0\degree]^T \\
\textbf{\oe}_2&=[4/1,0,50\degree,0\degree,170.2\degree,0\degree]^T
\end{split}
\end{equation*}

In this case, the first satellite of the system is essentially identical to the seed satellite (i.e., a unit impulse at $n=0$). The circular convolution between the seed satellite access profile and the constellation pattern vector yields the coverage timeline shown in the bottom part of Fig.~\ref{fig:example:a}. A snapshot of the corresponding configuration in the Earth-centered inertial (ECI) and Earth-centered Earth-fixed (ECEF) frame at $n=0$ is shown in Fig.~\ref{fig:example:b}.

This formulation exhibits the satellite constellation architecture by laying out the relationships between the common orbital characteristics, the satellite constellation pattern, and the coverage performance. We shall hereafter refer to this type of satellite constellation design decomposition into three vectors $\bm{v}_{0,j}$, $\bm{x}$, and $\bm{b}_j$ as the \textit{APC decomposition}, following the acronyms of the seed satellite Access profile, constellation Pattern, and Coverage timeline. Methods that are derived based on the APC decomposition are called the \textit{APC-based methods}.

\section{Regional Coverage Constellation Pattern Design Methods} \label{sec:methodology}
\subsection{Problem Statement}
Following the APC decomposition introduced herein, the satellite constellation design can be split into defining the reference seed satellite orbital elements $\textbf{\oe}_0$ (which includes the common orbital characteristics) and defining the constellation pattern vector $\bm{x}$. Conventional methods often make simple assumptions for $\bm{x}$ such as a symmetric pattern (e.g., Walker constellations) and optimize $\textbf{\oe}_0$; instead, this paper focuses on the optimization of the $\bm{x}$ itself without such simplifying assumptions. Mathematically, the goal of this paper is to solve for the optimal constellation pattern vector $\bm{x}^\ast$ such that the coverage timeline $\bm{b}_j^\ast=\bm{v}_{0,j}\circledast \bm{x}^\ast$ is equal to or greater than the designated $\bm{f}$ coverage threshold. The objective function is the number of satellites required $N$, which can be deduced from Eq.~\eqref{eq:n_sat}. The seed satellite orbital elements $\textbf{\oe}_0$ is considered as a given input so that the developed constellation pattern design approach can be integrated with the existing established methods for determining $\textbf{\oe}_0$ (e.g., brute-force methods, genetic algorithms). Appendix~C introduces an example approach to integrate the determination of the seed satellite orbital elements $\textbf{\oe}_0$ and the design of the satellite constellation pattern design $\bm{x}$.

This section introduces two constellation pattern optimization methods based on the circular convolution formulation and APC decomposition. First, we derive a rather conventional iterative method using a common assumption of symmetry; this method is used as a baseline for later analysis. Next, we develop a novel and general method based on binary integer linear programming to perform rigorous optimization of the constellation pattern.

\subsection{Baseline: Quasi-Symmetric Method}
The baseline quasi-symmetric method aims to design the satellite constellation pattern with uniform temporal spacing between satellites along the common closed trajectory in space. Given a length $L$ of the constellation pattern vector, the uniform temporal spacing constant $\eta\in \mathbb{R}_{>0}$ between satellites is defined as:
\begin{equation}
\label{eq:eta}
\eta\triangleq\frac{L}{N}
\end{equation}

We first consider a special case where $\eta$ is an integer. In this case and assuming $n_1=0$, we can construct a symmetric constellation pattern (i.e., a uniform distribution of satellites along the common ground track of the constellation system) using the following the constellation pattern vector form:

\begin{equation}
\label{eq:x_bar}
\bar{x}[n]\triangleq\sum_{k=1}^{N}\delta[n-\eta(k-1)]
\end{equation}
where 
\begin{equation}
    \delta[n] = \begin{cases}
	1, &\text{if } n=0 \\
	0, &\text{otherwise}
	\end{cases}
\end{equation}
A user is allowed to arbitrarily set the temporal location of the first satellite $n_1$ ($0\le n_1<L$). In this case, Eq.~\eqref{eq:x_bar} requires a circular shift of $\bar{x}[n]$:

\begin{equation}
\label{eq:x_qs}
x[n]=\bar{x}[n] \circledast \delta[n-n_1]
\end{equation}

Next, we generalize this formulation into the case where $\eta$ is not an integer. In this case, we cannot achieve a strictly symmetric constellation pattern with the given discretization, but only a near-symmetric one; we call the latter a \textit{quasi-symmetric} constellation pattern in this paper. For this generalization, the only change we need to make is to replace Eq.~\eqref{eq:x_bar} by Eq.~\eqref{eq:x_barnint}:
\begin{equation}
\label{eq:x_barnint}
\bar{x}[n]\triangleq\sum_{k=1}^{N}\delta[\text{nint}(n-\eta(k-1))]
\end{equation}
where nint($\cdot$) is the nearest integer function, which is used to guarantee the integer-indexing of a vector.

Algorithm~\ref{alg:qs} is designed to perform an iterative search about $N$ and $n_1$ until the coverage requirement is satisfied and outputs the optimal constellation pattern vector $\bm{x}^\ast$ given a set of $\bm{v}_{0,j}$ and $\bm{f}_j$. The algorithm consists of two nested iterative loops. The outer loop increments $N$ by one at each iteration, whereas the inner loop performs an exhaustive search about $n_1$ to find the $N$-minimizing temporal location of the first satellite. (Note that the range for the inner loop is set to $0\leq n_1 \leq \text{nint}({\eta})-1$ due to the (quasi-)symmetry of the resulting constellation pattern vector.) These loops break when the coverage requirement is satisfied as shown in Algorithm~\ref{alg:qs}. If no quasi-symmetric constellation is found until the outer loop for $N$ reaches the maximum number of satellites, which is equal to $L$, the method would determine the problem to be infeasible. 

\begin{algorithm}[H]
	\caption{The quasi-symmetric method to compute $\bm{x}^\ast$, $\bm{b}_j^\ast$, $N$, and $n_1$ (point-coverage)}
	\begin{algorithmic}[1]
		\Procedure{Quasi-SymmetricMethod}{$\bm{v}_{0,j},\bm{f}_j$}
		\State $N=1$
		\While {True}
		\If {$N\le L$}
    	\State Generate $\bar{x}[n]$ based on $\eta\triangleq L/N$ as outlined in Eq.~\eqref{eq:x_bar}
    	\For {$n_1=0,...,\text{nint}({\eta})-1$}
    	\State Generate $x[n]$ based on $\bar{x}[n]$ and $n_1$ as outlined in Eq.~\eqref{eq:x_qs}
    	\State Compute $b_j[n] = v_{0,j}[n] \circledast x[n]$ via Eq.~\eqref{eq:circular}
    	\If {$c_j=1$ as in Eq.~\eqref{eq:c_j}}
    	\State Break the loops
    	\State \Return $\bm{x}^\ast$, $\bm{b}_j^\ast$, $N$ (Eq.~\eqref{eq:n_sat}), and $n_1$
    	\EndIf
    	\EndFor
    	\State $N=N+1$
    	\Else
    	\State \Return Infeasible
    	\EndIf
    	\EndWhile
		\EndProcedure
	\end{algorithmic}
\label{alg:qs}
\end{algorithm}

For an area of interest consisting of multiple target points, a user may replace line~9 in Algorithm~\ref{alg:qs} with ``\textbf{if} $c_{\mathcal{J}}=1$ as in Eq.~\eqref{eq:c} \textbf{then}''. This guarantees the iterative search until all target points are satisfactorily covered. Similarly, one can come up with a custom termination criterion and/or figure of merit, such as time percent coverage or area percent coverage metrics. This is feasible since each iteration provides a full coverage state across all target points.

\begin{figure}[H]
	\centering\includegraphics[width=0.6\textwidth]{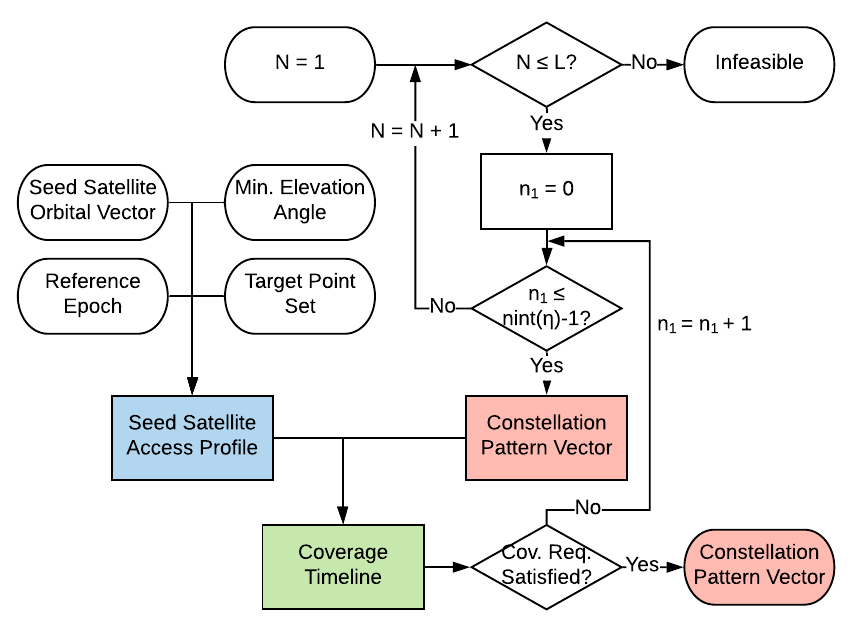}
	\caption{Overview of the quasi-symmetric method}
	\label{fig:qs}
\end{figure}

An overview of the quasi-symmetric method is shown in Fig.~\ref{fig:qs}. The seed satellite orbital elements vector, minimum elevation angle, reference epoch, and a set of target points are the user-defined parameters, which are determined based on mission requirements.

\subsection{New Method: Binary Integer Linear Programming (BILP) Method} \label{sec:bilp}
This subsection introduces the new satellite constellation pattern method developed in this paper using BILP. The BILP method aims to optimize the constellation pattern in a more rigorous and general way, without assuming symmetry and, if needed, concurrently considering multiple sub-constellations. Recall Eq.~\eqref{eq:circulantX}:

\begin{equation}
\label{eq:circulant}
\bm{V}_{0,j}\bm{x}=\bm{b}_j
\end{equation}

{\parindent0pt
	where $\bm{V}_{0,j}\in\mathbb{Z}^{L\times L}_{2}$ is a circulant matrix that is fully specified by the seed satellite access profile $\bm{v}_{0,j}$ as shown in Eq.~\eqref{eq:vdef}. This definition of $\bm{V}_{0,j}$ can be expanded as Eq.~\eqref{eq:V}.
}

\begin{equation}
	\label{eq:V}
	\bm{V}_{0,j}=
	\begin{bmatrix}
		v_{0,j}[0] & v_{0,j}[L-1] & v_{0,j}[L-2] & \cdots &v_{0,j}[1] \\
		v_{0,j}[1] & v_{0,j}[0] & v_{0,j}[L-1] & \cdots & v_{0,j}[2] \\
		v_{0,j}[2] & v_{0,j}[1] & v_{0,j}[0] & & \vdots \\
		\vdots & \vdots & \ddots & \ddots & \\
		v_{0,j}[L-1] & v_{0,j}[L-2] & \cdots & & v_{0,j}[0]
	\end{bmatrix}
\end{equation}
Each column of a circulant matrix $\bm{V}_{0,j}$ is identical to a circularly-shifted seed satellite access profile $\bm{v}_{0,j}$. Eq.~\eqref{eq:circulant} can be shown in a matrix form,

\begin{equation}
	\renewcommand{\arraystretch}{1.2}
	\begin{bmatrix}
		v_{0,j}[0] & v_{0,j}[L-1] & v_{0,j}[L-2] & \cdots &v_{0,j}[1] \\
		v_{0,j}[1] & v_{0,j}[0] & v_{0,j}[L-1] & \cdots & v_{0,j}[2] \\
		\vdots & \vdots & \ddots & \ddots & \vdots \\
		v_{0,j}[L-1] & v_{0,j}[L-2] & \cdots & & v_{0,j}[0]
	\end{bmatrix}
	\begin{bmatrix}
		x[0] \\
		x[1] \\
		\vdots \\
		x[L-1]
	\end{bmatrix}
	=
	\begin{bmatrix}
		b_j[0] \\
		b_j[1] \\
		\vdots \\
		b_j[L-1]
	\end{bmatrix}
\end{equation}

An interesting observation can be formalized. If we are given $\bm{v}_{0,j}$ and $\bm{x}$, then we can produce $\bm{b}_j$---this is the assumption of the quasi-symmetric method at each iteration. Likewise, if $\bm{v}_{0,j}$ and $\bm{b}_j$ are given, then we can analytically solve for $\bm{x}$ by solving the system of linear equations in Eq.~\eqref{eq:circulant} and obtain $\bm{x}=\bm{V}_{0,j}^{-1}\bm{b}_j$ ($\det(\bm{V}_{0,j})\neq0$). Since $\bm{b}_j$ represents the entire coverage timeline, this analysis enables us to find a constellation pattern vector $\bm{x}$ that satisfies a given coverage requirement $\bm{f}_j$.

Although this approach provides us with a way to find the satellite constellation pattern, the resulting $\bm{x}$ is not necessarily a binary vector, which violates the nature of the constellation pattern vector. The existence of a satellite at a given instance cannot be represented in a decimal number but only as either one or zero. Therefore, to guarantee a physical quantification of satellites, we shall employ the binary integer linear programming, or BILP, to solve for $\bm{x}^\ast$ which satisfies the inequality constraint:

\begin{equation}
\label{eq:inequality}
\bm{V}_{0,j}\bm{x}^\ast=\bm{b}_j^\ast \ge \bm{f}_j
\end{equation}

Before we formalize the BILP problem that solves Eq.~\eqref{eq:inequality}, Sections~\ref{sec:mtp} and \ref{sec:msc} introduce linear properties associated with Eq.~\eqref{eq:circulant}.

\subsubsection{Multiple Target Points} \label{sec:mtp}

Because the system is linear, we can extend Eq.~\eqref{eq:circulant} to an area of interest that consists of \textit{multiple target points}.

\begin{equation}
\label{eq:augmented}
\renewcommand{\arraystretch}{1}
\begin{bmatrix}
\bm{V}_{0,1} \\
\bm{V}_{0,2} \\
\vdots \\
\bm{V}_{0,\left\lvert \mathcal{J} \right\rvert} \\
\end{bmatrix}
\bm{x}
=
\begin{bmatrix}
\bm{b}_1 \\
\bm{b}_2 \\
\vdots \\
\bm{b}_{\left\lvert \mathcal{J} \right\rvert} \\
\end{bmatrix}
\end{equation}

{\parindent0pt
	where $\left\lvert \mathcal{J} \right\rvert$ is the cardinality of a target point set $\mathcal{J}$. Eq.~\eqref{eq:augmented} has the dimension of $({\left\lvert \mathcal{J} \right\rvert} L \times L)\cdot(L \times 1)=({\left\lvert \mathcal{J} \right\rvert} L \times 1)$.
}

The augmented circulant matrix on the left-hand side is a matrix of matrices obtained by appending all circulant matrices $\bm{V}_{0,1},...,\bm{V}_{0,\left\lvert \mathcal{J} \right\rvert}$ linearly. Similarly, the augmented coverage timeline vector is also obtained by appending all coverage timeline vectors $\bm{b}_1,...,\bm{b}_{\left\lvert \mathcal{J} \right\rvert}$ linearly. Here, the constellation pattern vector $\bm{x}$ represents a single constellation configuration that satisfies the augmented linear condition.

\subsubsection{Multiple Sub-Constellations} \label{sec:msc}
Another direction of linearity regarding having multiple sub-constellations is observed. We consider a constellation system consisting of \textit{multiple sub-constellations} with different seed satellite access profiles, $\bm{v}_{0,j}^{(1)},...,\bm{v}_{0,j}^{(z)},...,\bm{v}_{0,j}^{(\left\lvert \mathcal{Z} \right\rvert)}$, where superscript $z$ in parenthesis denotes the index of a sub-constellation, $\mathcal{Z}$ is a set of sub-constellations, and $\left\lvert \mathcal{Z} \right\rvert$ represents its cardinality. Each sub-constellation seed satellite access profile $\bm{v}_{0,j}^{(z)}$ is computed based on its seed satellite orbital elements vector $\textbf{\oe}_0^{(z)}$ and the modified minimum elevation angle threshold $\varepsilon_{j,\text{min}}^{(z)}[n], \ j \in \mathcal{J}, \ z \in \mathcal{Z}, \ n\in\{0,...,L-1\}$, which is only applicable to the BILP method (since the quasi-symmetric method does not define multiple sub-constellations). The goal of the multiple sub-constellation system is to satisfy a common coverage requirement over a single target point $j$. Thus, this can be incorporated by replacing Eq.~\eqref{eq:circulant} by the following equation:

\begin{equation}
\label{eq:multisubcon}
\renewcommand{\arraystretch}{1}
\begin{bmatrix}
\bm{V}_{0,j}^{(1)} & \bm{V}_{0,j}^{(2)} & \cdots & \bm{V}_{0,j}^{(\left\lvert \mathcal{Z} \right\rvert)}
\end{bmatrix}
\begin{bmatrix}
\bm{x}^{(1)} \\
\bm{x}^{(2)} \\
\vdots \\
\bm{x}^{(\left\lvert \mathcal{Z} \right\rvert)}
\end{bmatrix}
=
\bm{b}_j
\end{equation}

{\parindent0pt
	where the dimension of the system is $(L \times \left\lvert \mathcal{Z} \right\rvert L) \cdot (\left\lvert \mathcal{Z} \right\rvert L \times 1) = (L \times 1)$.
}

To guarantee the validity of this approach, we assume a synchronization condition among the sub-constellations to guarantee synchronized repeatability of the resulting coverage timeline:
\begin{equation}
T_{\text{r}}^{(1)}=...=T_{\text{r}}^{(z)}=...=T_{\text{r}}^{(\left\lvert \mathcal{Z} \right\rvert)}
\end{equation}
where $T_{\text{r}}^{(z)}, \ z\in \mathcal{Z}$ is the period of repetition, which can be written as a function of $a$, $e$, and $i$ and is therefore unique to each sub-constellation. Note that this does not mean that the individual orbital elements for each sub-constellation need to be all identical; instead, it only means that the period of repetition, defined by Eq.~\eqref{eq:period_of_repetition}, needs to be identical.

\subsubsection{A System of Multiple Sub-Constellations for Multiple Target Points}

Combining both directions of linearity---multiple target points and multiple sub-constellations---we get the following generalized governing relationship:

\begin{equation}
\label{eq:generalized}
\renewcommand{\arraystretch}{1.2}
\begin{bmatrix}
\bm{V}_{0,1}^{(1)} & \bm{V}_{0,1}^{(2)} & \cdots & \bm{V}_{0,1}^{(\left\lvert \mathcal{Z} \right\rvert)} \\
\bm{V}_{0,2}^{(1)} & \bm{V}_{0,2}^{(2)} & \cdots & \bm{V}_{0,2}^{(\left\lvert \mathcal{Z} \right\rvert)} \\
\vdots & \vdots & \ddots & \vdots \\
\bm{V}_{0,\left\lvert \mathcal{J} \right\rvert}^{(1)} & \bm{V}_{0,\left\lvert \mathcal{J} \right\rvert}^{(2)} & \cdots & \bm{V}_{0,\left\lvert \mathcal{J} \right\rvert}^{(\left\lvert \mathcal{Z} \right\rvert)} \\
\end{bmatrix}
\begin{bmatrix}
\bm{x}^{(1)} \\
\bm{x}^{(2)} \\
\vdots \\
\bm{x}^{(\left\lvert \mathcal{Z} \right\rvert)}
\end{bmatrix}
=
\begin{bmatrix}
\bm{b}_1 \\
\bm{b}_2 \\
\vdots \\
\bm{b}_{\left\lvert \mathcal{J} \right\rvert}
\end{bmatrix}
\end{equation}

{\parindent0pt
	where the dimension of the system is $(\left\lvert \mathcal{J} \right\rvert L\times \left\lvert \mathcal{Z} \right\rvert L) \cdot (\left\lvert \mathcal{Z} \right\rvert L \times 1) = (\left\lvert \mathcal{J} \right\rvert L\times 1)$.
}

Eq.~\eqref{eq:generalized} can be expressed in an indexed equation form:

\begin{equation}
\label{eq:indexed}
\sum_{z=1}^{\left\lvert \mathcal{Z} \right\rvert}\bm{V}_{0,j}^{(z)} \bm{x}^{(z)}=\bm{b}_j, \ \ \forall j\in\mathcal{J}
\end{equation}

{\parindent0pt
	where the subscript $j$ is the target point index and the superscript $z$ is the sub-constellation index.
}

The physical interpretation of Eq.~\eqref{eq:generalized} is as follows: it represents a linear relationship between the physical configuration of a system of multiple sub-constellations and the resulting coverage timelines over a set of multiple target points. Here, each sub-constellation may exhibit its own unique orbital characteristics. For example, a sub-constellation ($z=1$) may be placed on a critically-inclined elliptic orbit while a sub-constellation ($z=2$) may be placed on a circular low Earth orbit. Similarly, each target point may impose an independent coverage requirement. For example, a target point ($j=1$) may require continuous single-fold coverage whereas a target point ($j=2$) may require a sinusoidal-like time-varying coverage, fluctuating between the double and triple folds. Revisiting the inequality constraint as shown in Eq.~\eqref{eq:inequality}, it is the goal of the binary integer linear programming to determine the satellite constellation configurations $\bm{x}^{(1)},...,\bm{x}^{(z)},...,\bm{x}^{(\left\lvert \mathcal{Z} \right\rvert)}$ that satisfy this complex relationship.

\subsubsection{Binary Integer Linear Programming (BILP) Problem Formulation}
Let us assume that we want to achieve a $\bm{f}_j$-fold coverage system ($\forall j\in\mathcal{J}$) with the given $\bm{v}_{0,j}^{(1)},...,\bm{v}_{0,j}^{(z)},...,\bm{v}_{0,j}^{(\left\lvert \mathcal{Z} \right\rvert)}$ vectors. The BILP formulation is shown in Eq.~\eqref{eq:bilp}. Solving the most general form of the problem, Eq.~\eqref{eq:indexed}, via BILP yields an optimal solution in the form of ``\textit{a system of multiple sub-constellations that simultaneously satisfies the coverage requirements over multiple target points}.''

\begin{equation}
\label{eq:bilp}
\begin{aligned}
& \underset{\bm{x}}{\text{minimize}}
& & \bm{1}^T\bm{x} \\
& \text{subject to}
& & \sum_{z=1}^{\left\lvert \mathcal{Z} \right\rvert}\bm{V}_{0,j}^{(z)} \bm{x}^{(z)} \ge \bm{f}_j, \ \ \forall j\in\mathcal{J} \\
&
& & \bm{x}\in\mathbb{Z}^{L}_{2} \\
\end{aligned}
\end{equation}

{\parindent0pt
	where the binary design variable constraint is imposed on the elements of the constellation pattern vector $\bm{x}$ to reflect the physical quantification of satellites. The solution to this BILP problem is the optimal constellation pattern vector $\bm{x}^\ast$.
}

\begin{figure}[H]
	\centering\includegraphics[width=0.55\textwidth]{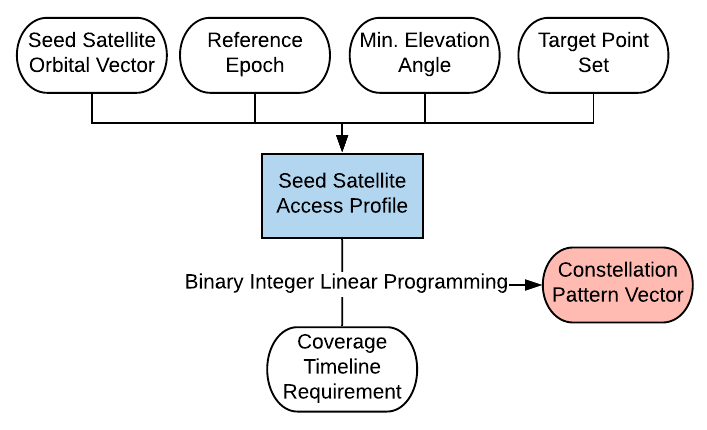}
	\caption{Overview of the binary integer linear programming method}
	\label{fig:ip}
\end{figure}

An overview of the BILP method is shown in Fig.~\ref{fig:ip}.

\subsection{Derivation of $\Omega$ and $M$ from the Constellation Pattern Vector} \label{sec:post}

Once the aforementioned methods obtain an optimal constellation pattern vector, it must be post-processed to extract interpretable orbital information---a set of $(\Omega,M)_k$ where $k$ is the index of a satellite. Every impulse on a constellation pattern vector corresponds to a point in the $(\Omega,M)$-space. Given $n_k$ found from the constellation pattern vector, one can find $(\Omega,M)_k$ set by solving the following system of equations:
\begin{subequations}
\label{eq:set}
\begin{align}
N_{\text{P}} (\Omega_k-\Omega_0)+N_{\text{D}} (M_k-M_0) = 0 \ \text{mod} \ (2\pi) \label{eq:seta} \\
\Omega_k=n_k\frac{2\pi N_{\text{D}}}{L}+\Omega_0 \label{eq:setb}
\end{align}
\end{subequations}
Note that Eq.~\eqref{eq:seta} is rearranged from Eq.~\eqref{eq:flower} \cite{avendano2013}. The derivation of Eq.~\eqref{eq:setb} is explained in Appendix~D.

\section{Illustrative Examples} \label{sec:illustrative_examples}
This section aims to demonstrate the general applicability of and the computational efficiency associated with the proposed methods under various mission profiles. Five illustrative examples are uniquely set up by varying orbital characteristics, area of interest properties, minimum elevation angle, and coverage requirements to illustrate the APC decomposition.

All illustrative examples are conducted on an Intel Core i9-9940X Processor @3.30 GHz platform. For BILP problems, Gurobi 9.0.0 is used with the default termination setting \cite{gurobi}. The referenced ellipsoid model adopts the World Geodetic System 1984 (WGS 84). 
It is assumed that all satellites point to their nadir directions. 
Furthermore, we assume the utilization of satellite maneuvers to correct and maintain an identical ground track throughout the satellite lifetime, negating the perturbation effects other than the $J_2$ effect. Lastly, we make an assumption that the minimum elevation angle threshold is time-invariant:

\begin{equation*}
    \varepsilon_{j,\text{min}}[n]=\varepsilon_{j,\text{min}}, \ \ \ \forall n \in\{0,...,L-1\}
\end{equation*}

Table~\ref{tab:examples} is a list of parameters used for each example study. The five examples are chosen to test different capabilities of the methods: Example 1 for single-fold continuous coverage over a single target point; Example 2 for time-varying coverage over a single target point; Example 3 for single-fold continuous coverage over multiple target points; Example 4 for time-varying and spatially-varying coverage over multiple target points; and Example 5 for multiple sub-constellations over multiple target points. All examples uniquely illustrate a variety of orbit (circular vs. critically-inclined elliptic, prograde vs. retrograde, and low vs. high altitudes) and a variety of areas of interest (a single target point vs. multiple target points and contiguous vs. discontiguous). Both the baseline quasi-symmetric method and the BILP method are applied to all examples, with an exception of the quasi-symmetric method for Example 5 due to its incapability of handling multiple sub-constellations. In this section, the subscripts \textit{qs} and \textit{bilp} denote variables associated with the quasi-symmetric and the BILP methods, respectively. The rest of this section discusses the details of each illustrative case.

\begin{table}[h]
		\fontsize{10}{10}\selectfont
		\caption{Example parameters}
		\renewcommand{\arraystretch}{2}
		\centering
		\begin{threeparttable}
		\begin{adjustbox}{max width=1.1\textwidth,center}
		\begin{tabular}{c c c c c c}
			\hline 
			\hline
			Example & Seed Satellite Orbital Elements \tnote{a,b} & Min. Elev. Angle & Target Point Set & Cov. Req. & L \\
			\hline
            1 & $[12/1,0,102.9\degree,0\degree,98.3\degree,0\degree]^T$ & $5\degree$ & $\{(\phi=34.75\degree \text{N},\lambda=84.39\degree \text{W})\}$ & $\bm{1}$ & 720 \\
			2 & $[12/1,0,102.9\degree,0\degree,98.3\degree,0\degree]^T$ & $5\degree$ & $\{(\phi=34.75\degree \text{N},\lambda=84.39\degree \text{W})\}$ & Time-Varying & 720 \\
			3 & $[5/1,0.41,63.435\degree,90\degree,0\degree,0\degree]^T$ & $30\degree$ & $\{\text{Antarctica}\}$ & $\bm{1}$ & 718 \\
			\multirow{2}{*}{4} & \multirow{2}{*}{$[83/6,0,99.2\degree,0\degree,0\degree,0\degree]^T$} & \multirow{2}{*}{$20\degree$} & $\mathcal{J}_1=\{\text{Amazon River Basin}\}$, & Time-Varying \& & \multirow{2}{*}{4200} \\
			& & & $\mathcal{J}_2=\{\text{Nile River Basin}\}$ & Spatially-Varying \\
			\multirow{2}{*}{5} & $\textbf{\oe}_0^{(1)}=[8/1,0,70\degree,0\degree,0\degree,0\degree]^T$ & $\varepsilon_{1,\text{min}}=15\degree$ & $\{(\phi=64.14\degree \text{N},\lambda=21.94\degree \text{W}),$ & \multirow{2}{*}{$\bm{1}$} & \multirow{2}{*}{717} \\
			& $\textbf{\oe}_0^{(2)}=[6/1,0,47.915\degree,0\degree,0\degree,0\degree]^T$ & $\varepsilon_{2,\text{min}}=10\degree$ & $(\phi=19.07\degree \text{N},\lambda=72.87\degree \text{E})\}$ \\
			\hline
			\hline
		\end{tabular}
		\end{adjustbox}
	\begin{tablenotes}
		\item[a] The seed satellite orbital elements vector $\textbf{\oe}_0$ takes the form of $[\tau,e,i,\omega,\Omega_0,M_0]^T$.
		\item[b] All orbital elements are in J2000.
		\end{tablenotes}
\end{threeparttable}
		\label{tab:examples}
\end{table}

\subsection{Example 1. Single-Fold Continuous Coverage over a Single Target Point}
A target point is located at $\{(\phi=34.75\degree \text{N},\lambda=84.39\degree \text{W})\}$ and requires $\varepsilon_{\text{min}}=5\degree$. A seed satellite orbital elements vector $\textbf{\oe}_0=[12/1,0,102.9\degree,0\degree,98.3\degree,0\degree]^T$ is assumed. The period of repetition is \SI{86400}{s}. The length of vectors is selected, $L=720$, such that the time step is $\SI{120}{s}$. The objective is to find the optimal constellation pattern vector $\bm{x}^\ast$ that satisfies a single-fold continuous coverage requirement ($\bm{f}=\bm{1}$).

The results are obtained as follows:
\begin{subequations}
	\begin{align*}
	x_{\text{qs}}^\ast[n]&=\begin{cases}
	1, &\text{for} \ n= 	\begin{aligned}
	&0,33,65,98,131,164,196,229,262,295,327,360,393,425,458,491,524,\ldots \\
	&556,589,622,655,687
	\end{aligned}\\
	0, &\text{otherwise}
	\end{cases} \\
	x_{\text{bilp}}^\ast[n]&=\begin{cases}
	1, &\text{for} \ n= 39,73,79,89,170,184,234,250,331,341,347,492,502,542,638,648,654,663 \\
	0, &\text{otherwise}
	\end{cases}
	\end{align*}
\end{subequations}
where the total number of satellites obtained for each method is $N_{\text{qs}}=22$ and $N_{\text{bilp}}=18$, with the computational time \SI{0.1}{s} for the quasi-symmetric method and \SI{5937.6}{s} for the BILP method.
The results indicate that, although the computational cost for the BILP method is longer, it can explore a substantially larger design space and achieve a fewer-satellite configuration than the quasi-symmetric method by breaking the symmetry. 

\begin{figure}[H]
	\centering\includegraphics[width=0.6\textwidth]{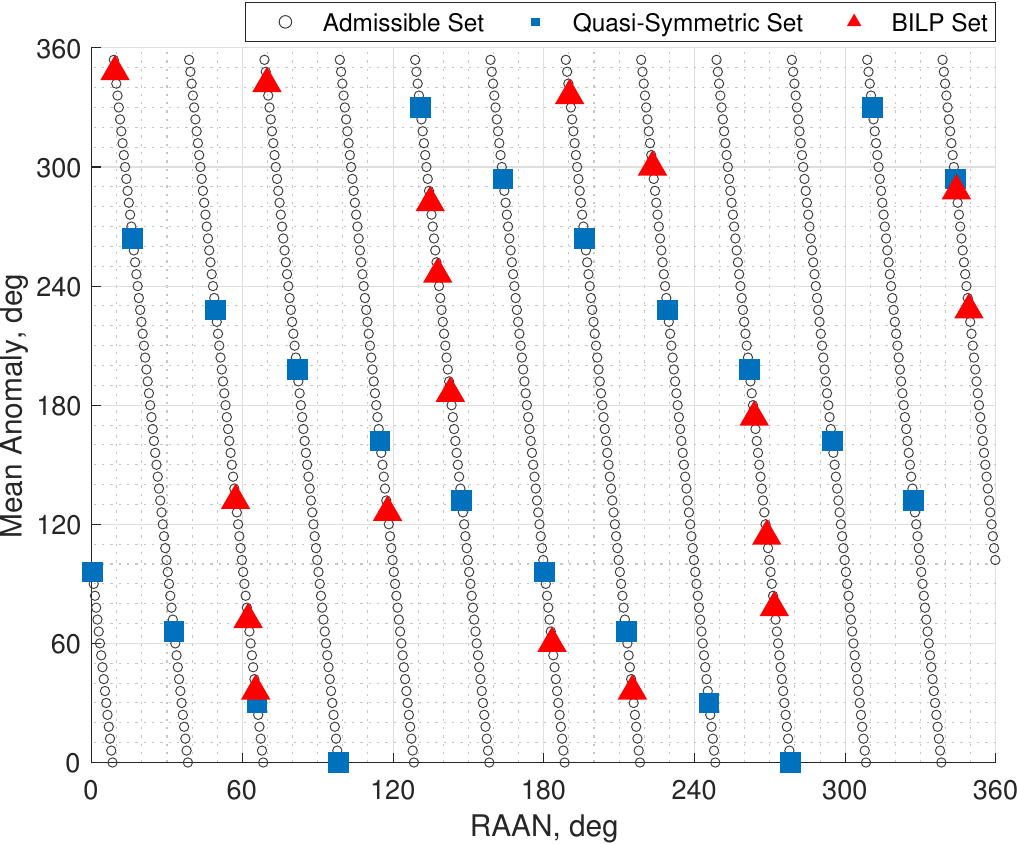}
	\caption{Example 1: Admissible set, quasi-symmetric set, and binary integer linear programming set in the $(\Omega,M)$-space}
	\label{fig:raan_m}
\end{figure}

Fig.~\ref{fig:raan_m} illustrates the $(\Omega,M)$-space and where each of the quasi-symmetric and BILP solution constellations lies. In this example, $L=720$; therefore, there are $L=720$ number of \textit{admissible} points in the $(\Omega,M)$-space into which a satellite can be placed. Analyzing the patterns in Fig.~\ref{fig:raan_m}, the quasi-symmetric set depicts a lattice-like symmetry in the $(\Omega,M)$-space whereas the BILP set exhibits asymmetry in the $(\Omega,M)$-space.

\begin{figure}[H]
	\centering
	\begin{subfigure}[h]{0.495\linewidth}
		\includegraphics[width=\linewidth]{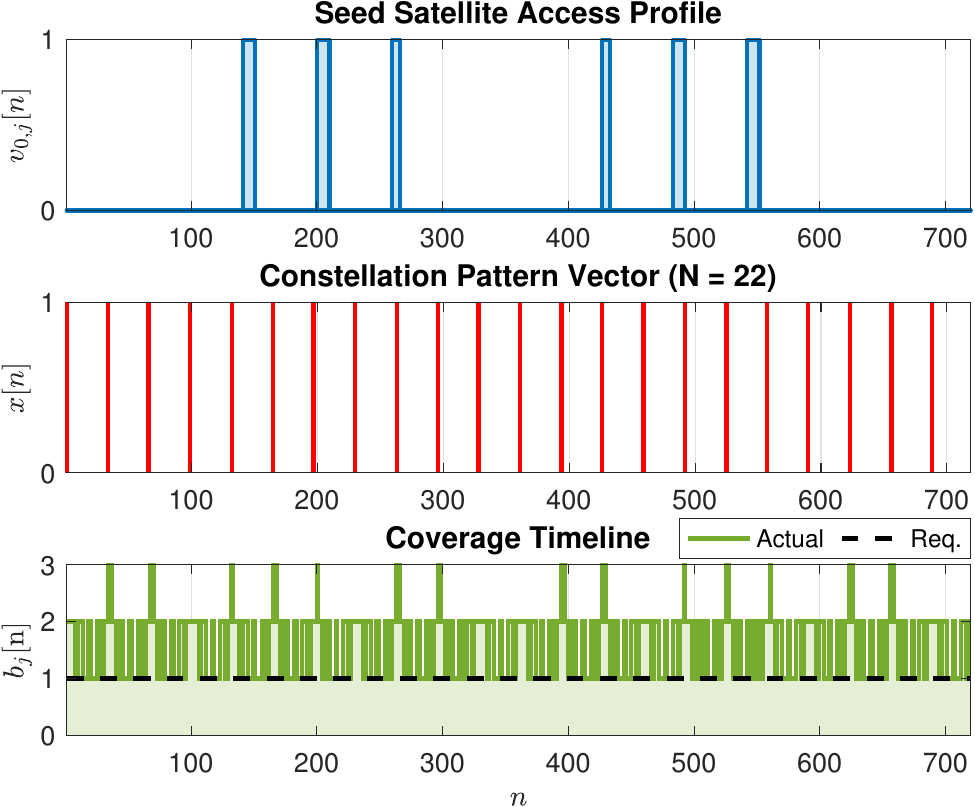}
		\caption{Quasi-Symmetric method}
	\end{subfigure}
	\begin{subfigure}[h]{0.495\linewidth}
		\includegraphics[width=\linewidth]{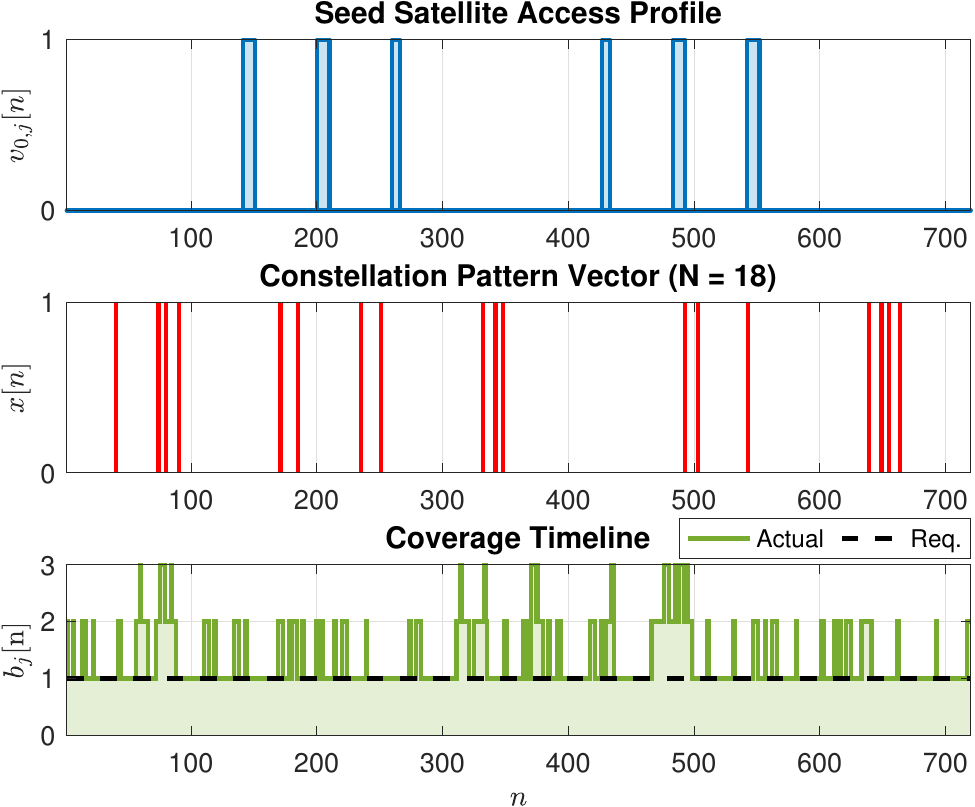}
		\caption{BILP method}
	\end{subfigure}
	\caption{Example 1: The APC decomposition}
	\label{fig:example1_apc}
\end{figure}

The APC decomposition figures are shown in Fig.~\ref{fig:example1_apc}. One can observe that the single-fold continuous coverage requirement is satisfied everywhere. Again, the asymmetry in the constellation pattern vector from the BILP method is contrasted with the symmetry in that from the quasi-symmetric method. Note that the coverage timeline for the quasi-symmetric constellation may not be strictly symmetric as $\eta$ is not an integer in this case.

\begin{figure}[H]
	\centering
	\begin{subfigure}[h]{0.4\linewidth}
	    \centering
		\includegraphics[height=5cm]{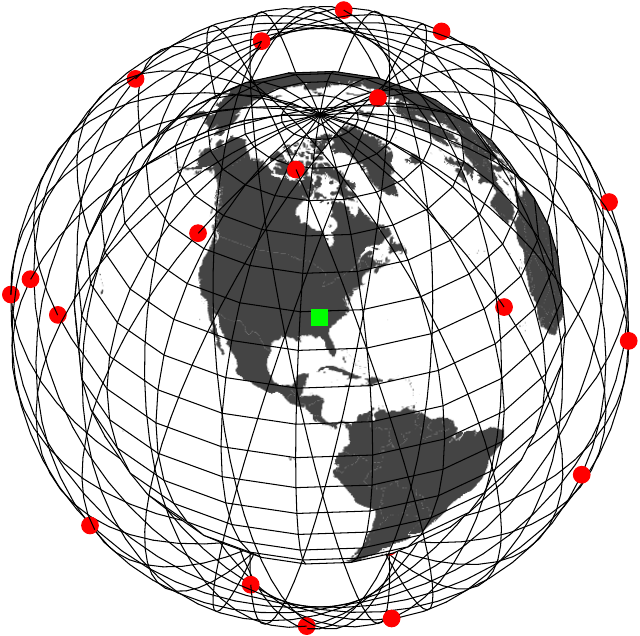}
		\caption{Quasi-Symmetric 22-sat constellation}
	\end{subfigure}
	\begin{subfigure}[h]{0.4\linewidth}
	    \centering
		\includegraphics[height=5cm]{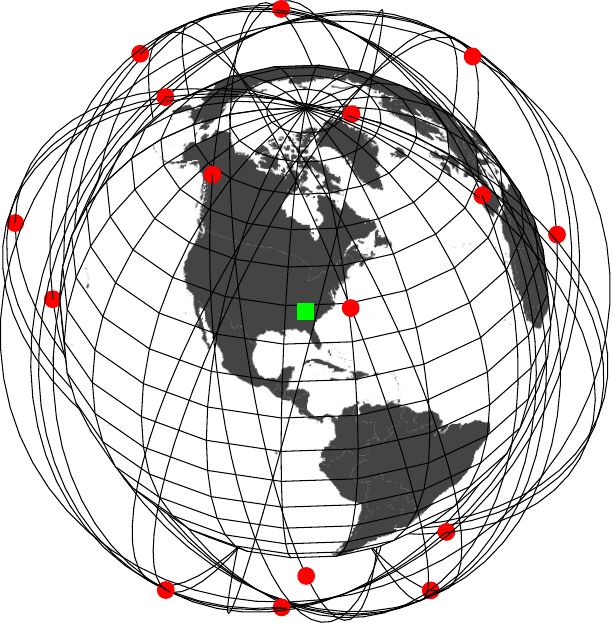}
		\caption{BILP 18-sat constellation}
	\end{subfigure}
	\caption{Example 1: 3D view of generated constellations at $n=0$ (ECI frame)}
	\label{fig:example1_3d}
\end{figure}

A snapshot of the corresponding constellation configurations at $n=0$ is shown in Fig.~\ref{fig:example1_3d}. This figure visually shows that the BILP method is taking advantage of the asymmetry to achieve a smaller number of satellites.

\subsection{Example 2. Time-Varying Coverage over a Single Target Point}
In this example, we execute a single variation to Example~1 such that the coverage requirement is now periodically time-varying with the rest of the parameters being identical (e.g., $L=720$). The objective is to find the optimal constellation pattern vector $\bm{x}^\ast$ that satisfies a specialized threshold function, namely, a square wave function:

\begin{equation*}
f[n]=\begin{cases}
2, & \text{for } 240 \le n \le 480 \\
1, &\text{otherwise}
\end{cases} \\
\end{equation*}

A coverage requirement is now time-dependent; the value of the square wave function varies between values 1 and 2. This requires that some parts of the simulation period must be continuously covered by at least two satellites (double-fold) and by at least one satellite (single-fold) during the other part of the simulation period. This case is an abstract illustration of general time-varying constellation applications. For example, a communication satellite constellation may require two satellites during the day for doubled-capacity and one satellite during the night for a quiescent mode.

The results are obtained as follows:
\begin{subequations}
	\begin{align*}
	x_{\text{qs}}^\ast[n]&=\begin{cases}
	1, &\text{for} \ n= 
	\begin{aligned}
	&0,22,44,65,87,109,131,153,175,196,218,240,262,284,305,327,349,\ldots \\
	&371,393,415,436,458,480,502,524,545,567,589,611,633,655,676,698
	\end{aligned}\\
	0, &\text{otherwise}
	\end{cases} \\
	x_{\text{bilp}}^\ast[n]&=\begin{cases}
	1, &\text{for} \ n=
	\begin{aligned}
	&5,23,39,75,89,114,124,130,164,215,230,255,265,483,493,518,533,\ldots \\
	&584,618,624,634,659,673,709
	\end{aligned}\\
	0, &\text{otherwise}
	\end{cases}
	\end{align*}
\end{subequations}

\begin{figure}[H]
	\centering\includegraphics[width=0.6\textwidth]{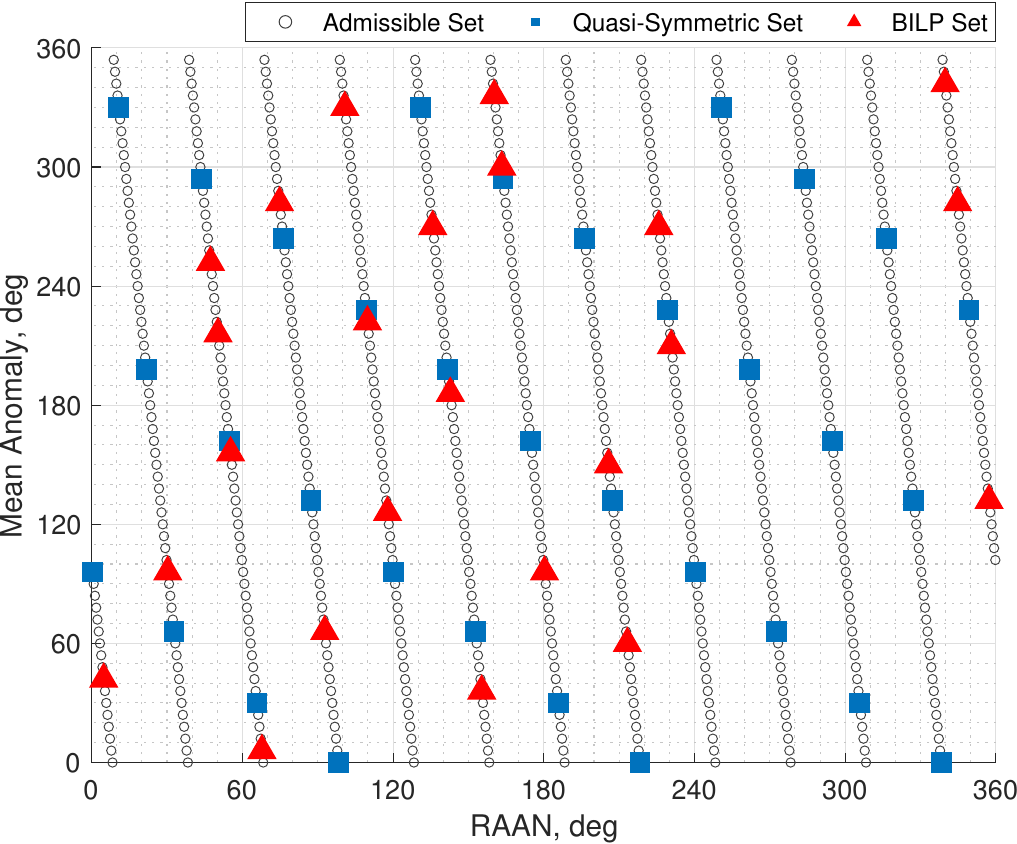}
	\caption{Example 2: Admissible set, quasi-symmetric set, and binary integer linear programming set in the $(\Omega,M)$-space}
	\label{fig:example2_raan_m}
\end{figure}

{\parindent0pt
	where the total number of satellites obtained for each method is $N_{\text{qs}}=33$ and $N_{\text{bilp}}=24$, and the computational time is \SI{0.1}{s} for the quasi-symmetric method and \SI{3712.0}{s} for the BILP method. Like in Example 1, although the BILP method takes longer computational time, it can achieve a constellation pattern solution that requires a significantly smaller number of satellites than the baseline quasi-symmetric method. The distribution of satellites in the $(\Omega,M)$-space is shown in Fig.~\ref{fig:example2_raan_m}.
}

\begin{figure}[H]
	\centering
	\begin{subfigure}[h]{0.495\linewidth}
		\includegraphics[width=\linewidth]{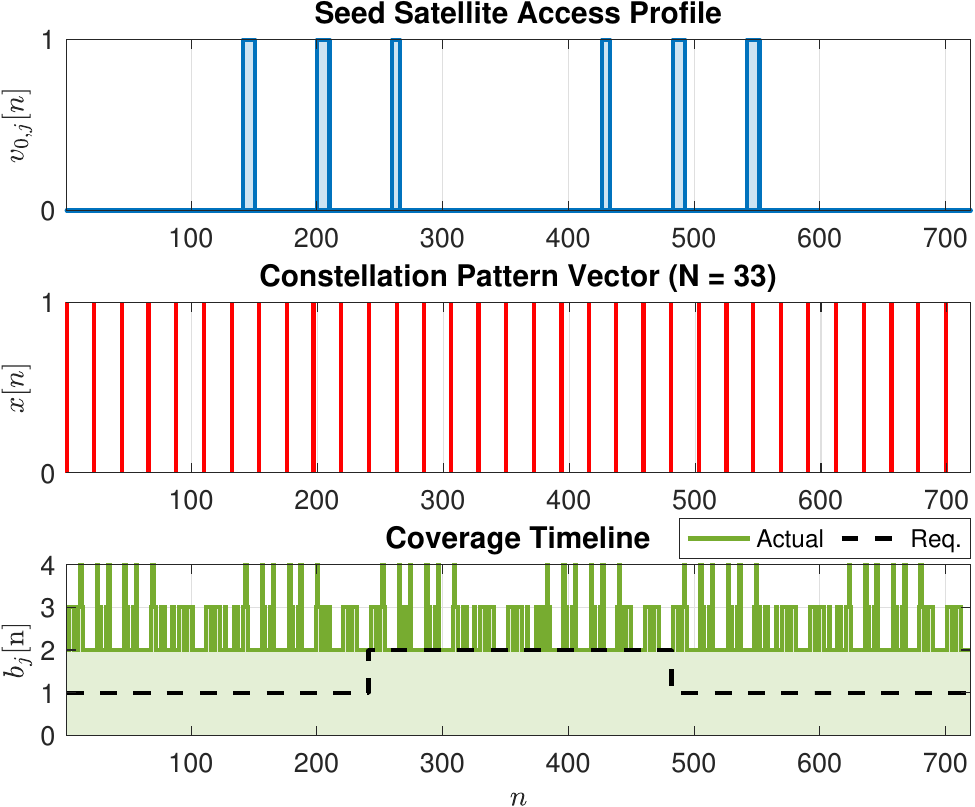}
		\caption{Quasi-Symmetric method}
	\end{subfigure}
	\begin{subfigure}[h]{0.495\linewidth}
		\includegraphics[width=\linewidth]{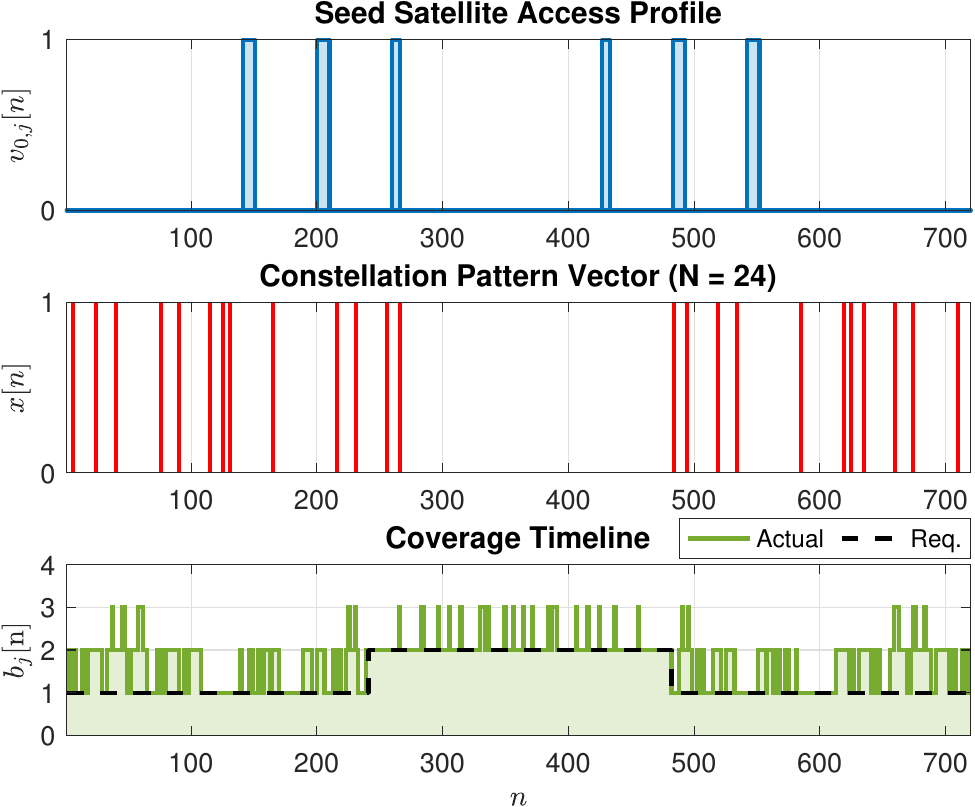}
		\caption{BILP method}
	\end{subfigure}
	\caption{Example 2: The APC decomposition}
	\label{fig:example2_apc}
\end{figure}

\begin{figure}[H]
	\centering
	\begin{subfigure}[h]{0.4\linewidth}
	    \centering
		\includegraphics[height=5cm]{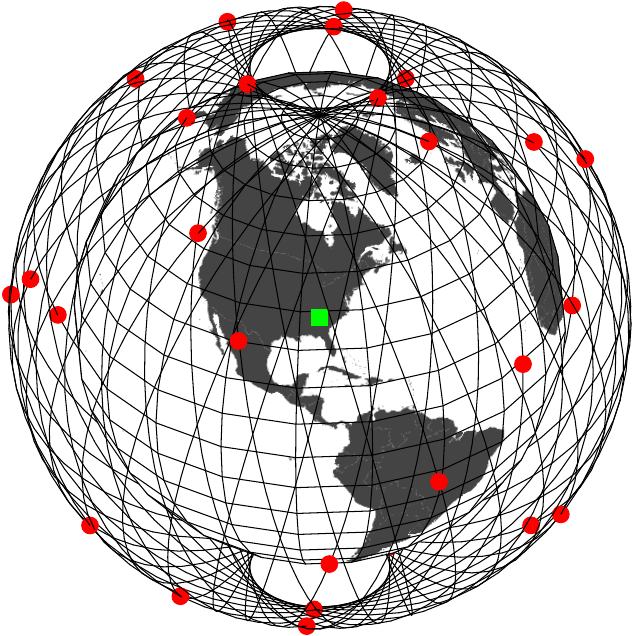}
		\caption{Quasi-Symmetric 33-sat constellation}
	\end{subfigure}
	\begin{subfigure}[h]{0.4\linewidth}
	    \centering
		\includegraphics[height=5cm]{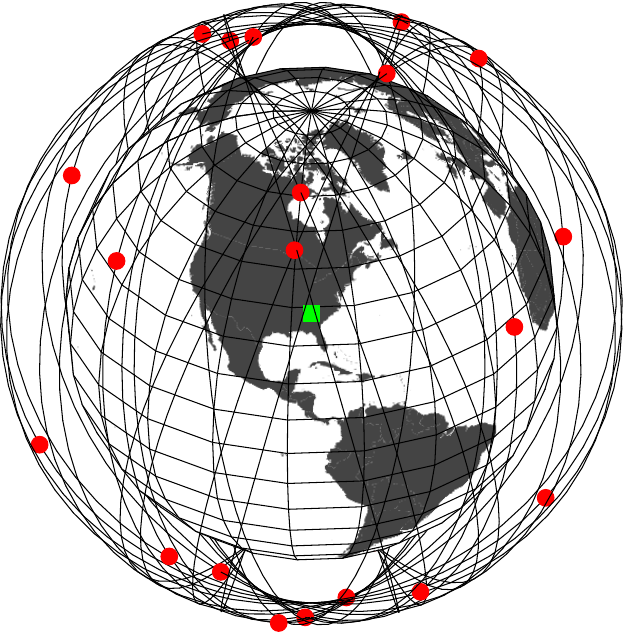}
		\caption{BILP 24-sat constellation}
	\end{subfigure}
	\caption{Example 2: 3D view of generated constellations at $n=0$ (ECI frame)}
	\label{fig:example2_3d}
\end{figure}

As shown in Fig.~\ref{fig:example2_apc}, the BILP constellation produces a coverage timeline that closely follows the time-varying coverage requirement. Such a coverage timeline is possible since the BILP constellation is not subject to symmetry in the satellite distribution. This is not the case for the quasi-symmetric method due to its (quasi-)symmetrical satellite distribution, which resulted in a conservative solution that provides a double-fold coverage over the entire period, even when it is not needed. This leads to the superior solution from the BILP method compared with the baseline quasi-symmetric method. As observed in Example 1, the BILP method already reduces the number of satellites required compared to that of the quasi-symmetric method given the single-fold coverage requirement. Changing only the coverage requirement to be time-varying, we further observe the additional reduction of the number of satellites for the BILP method. A snapshot of the corresponding constellation configurations at $n=0$ is shown in Fig.~\ref{fig:example2_3d}.

\subsection{Example 3. Single-Fold Continuous Coverage over Multiple Target Points}
For this example, we consider a target area, Antarctica, which calls for continuous and reliable telecommunication systems to support existing and planned scientific expeditions \cite{lee2016satellite}. The area is discretized into a set of 94 target points following the $3\degree$-by-$3\degree$ resolution (latitude-by-longitude). All target points set $\varepsilon_{\text{min}}=30\degree$. A seed satellite orbital element vector $\textbf{\oe}_0=[5/1,0.41,63.435\degree,90\degree,0\degree,0\degree]^T$ (critically-inclined elliptic orbit with the apogee over the southern hemisphere) is assumed. The period of repetition is \SI{86076}{s}. The length of vectors is selected, $L=718$, such that the time step is approximately $t_{\text{step}}\approx\SI{120}{s}$. The objective of this example is to design a satellite constellation configuration that achieves single-fold continuous coverage ($\bm{f}=\bm{1}$) over all target points.

\begin{figure}[H]
	\centering
	\includegraphics[trim=45 20 45 20,clip,width=0.55\linewidth]{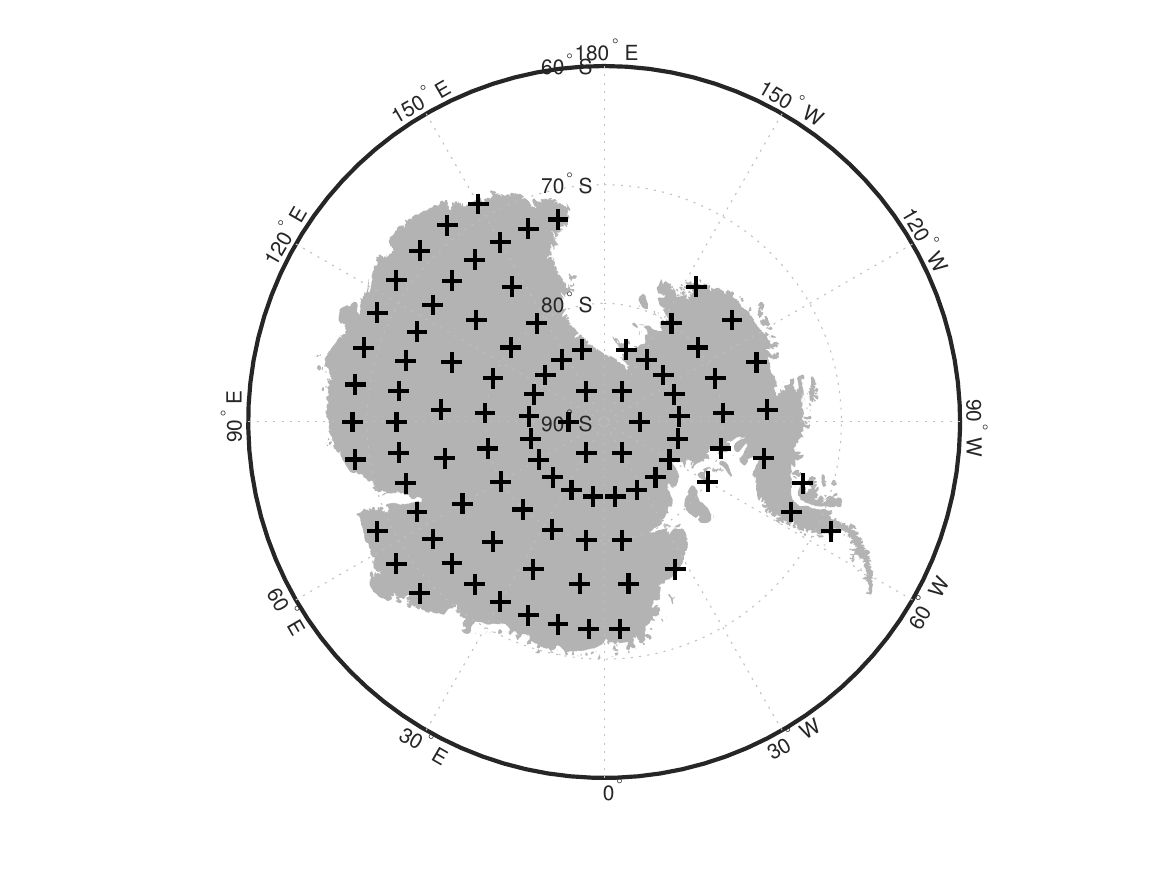}
	\caption{Example 3: Antarctica target points ($3\degree$-by-$3\degree$ resolution); the shapefile is obtained from Ref.~\cite{antarctica}}
	\label{fig:example3_area}
\end{figure}

Note that this continuous polar coverage is a typical example that is often handled with a symmetric constellation, and thus we would expect that the quasi-symmetric method would perform well. 

The results are obtained as follows:

\begin{subequations}
	\begin{align*}
	x_{\text{qs}}^\ast[n]&=\begin{cases}
	1, &\text{for} \ n=0,120,239,359,479,598 \\
	0, &\text{otherwise}
	\end{cases} \\
	x_{\text{bilp}}^\ast[n]&=\begin{cases}
	1, &\text{for} \ n=96,310,358,562,612 \\
	0, &\text{otherwise}
	\end{cases}
	\end{align*}
\end{subequations}

\begin{figure}[H]
	\centering
	\begin{subfigure}[h]{0.4\linewidth}
	    \centering
		\includegraphics[height=5cm]{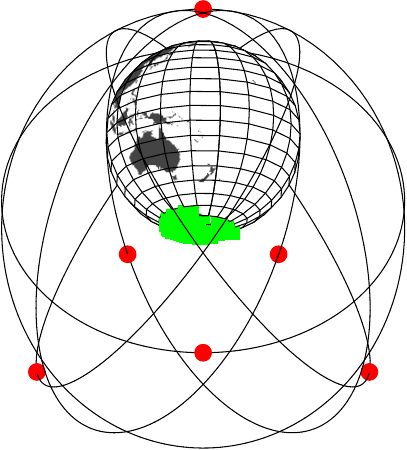}
		\caption{Quasi-Symmetric 6-sat constellation}
	\end{subfigure}
	\begin{subfigure}[h]{0.4\linewidth}
	    \centering
		\includegraphics[height=5cm]{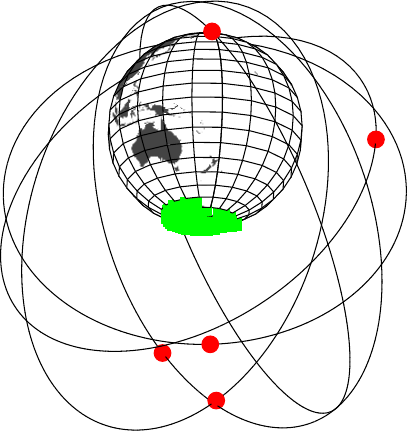}
		\caption{BILP 5-sat constellation}
	\end{subfigure}
	\caption{Example 3: 3D view of generated constellations at $n=0$ (ECI frame)}
	\label{fig:example3_3d}
\end{figure}

{\parindent0pt
	where the total number of satellites obtained for each method is $N_{\text{qs}}=6$ and $N_{\text{bilp}}=5$, and the computational cost was \SI{10.7}{s} for the quasi-symmetric method and \SI{748.4}{s} for the BILP method. It is worth mentioning that, even for this polar-coverage example for which we would typically just use a symmetric constellation pattern (i.e., using the baseline method), the BILP method still achieves an asymmetric constellation pattern with fewer satellites.  A snapshot of the obtained constellations is indicated in Fig.~\ref{fig:example3_3d}.
}

\subsection{Example 4. Time-Varying and Spatially-Varying Coverage over Multiple Target Points}
\begin{figure}[H]
	\centering\includegraphics[width=0.8\linewidth]{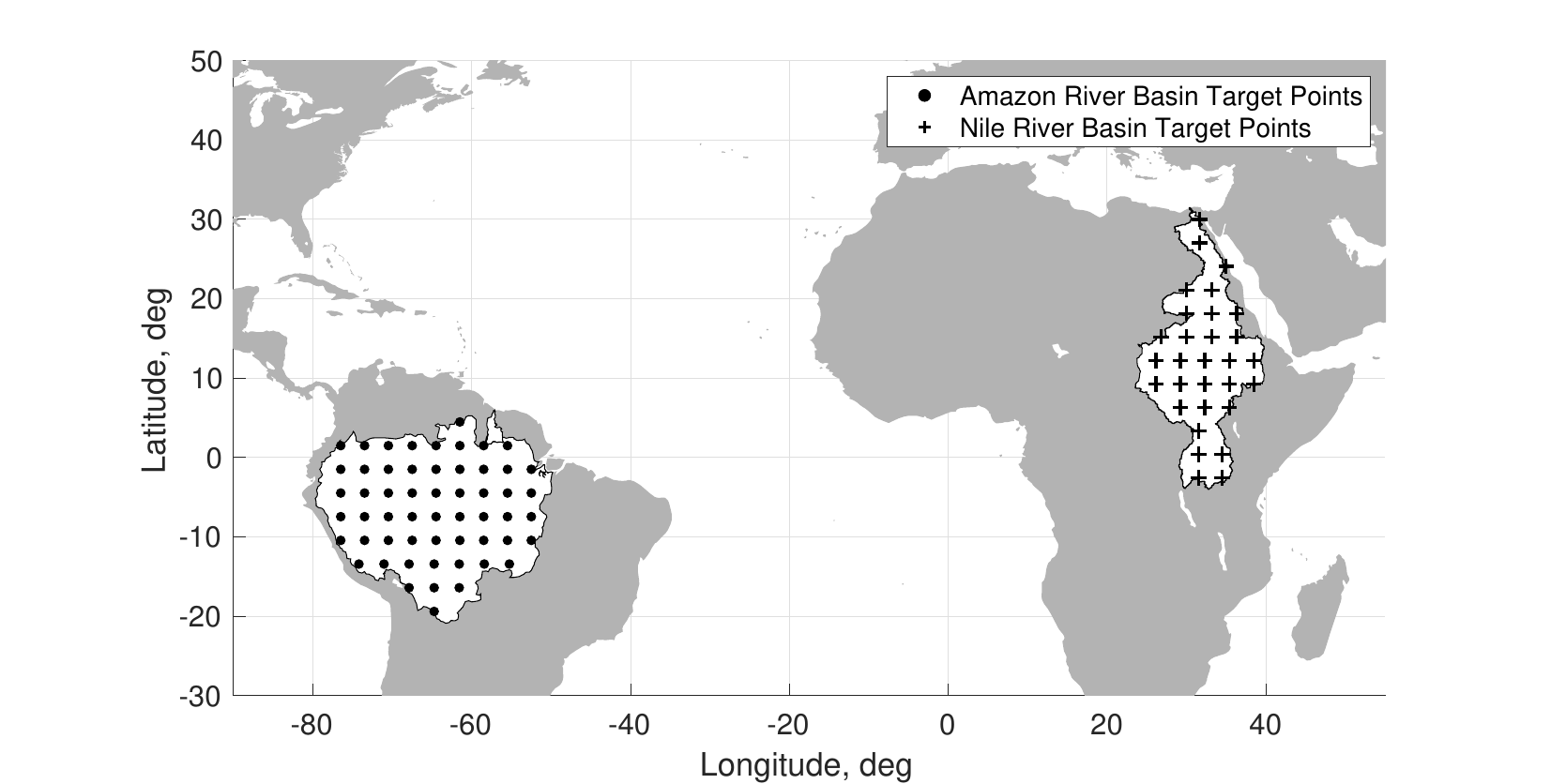}
	\caption{Example 4: Amazon and Nile river basin target points ($3\degree$-by-$3\degree$ resolution); the polygon shapefiles are retrieved from the dataset provided by the World Bank \cite{theworldbasin}}
	\label{fig:example4_area}
\end{figure}

In this example, we design a satellite constellation system that performs remote sensing tasks over two areas of interest: the Amazon and Nile river basins. These areas represent two of the major river basins in the world thereby making them desirable locations for monitoring forests, logging, soil and water managements \cite{chambers2007regional,rientjes2013diurnal}, and thus are of great interest to the international community. Each area of interest is discretized into a set of target points following the $3\degree$-by-$3\degree$ resolution (latitude-by-longitude). The Amazon river basin target point set $\mathcal{J}_1$ is composed of 56 target points and the Nile river basin target point set $\mathcal{J}_2$ is composed of 30 target points. The target points are shown in Fig.~\ref{fig:example4_area}.

Each target point set is assumed to require different revisit time requirements: the Amazon basin has a revisit time requirement of every twelve hours, starting six hours after the epoch, whereas the Nile basin has a revisit time requirement of every six hours, starting at the epoch. We assume that all target points within the same set require simultaneous access to the system satellites at given revisit time requirements. Note that these requirements are not just constraining the revisit time interval but the exact time step for revisit; this is referred to as the \textit{strict} revisit time requirement here. Furthermore, all target points are assumed to require the minimum elevation angle threshold of \SI{20}{\degree}, which corresponds to the hypothetical sensor's field-of-view of approximately \SI{110}{\degree} at a given altitude of satellites. The length of vectors is chosen, $L=4200$ ($t_{\text{step}}\approx\SI{123.4}{s}$), such that we can represent the complex coverage requirements in an integer-indexed symmetrical form. Note that this coverage requirement is both time-varying (i.e., periodic) and spatially-varying (i.e., different requirements for Amazon and Nile river basin target points).

\begin{subequations}
	\begin{align*}
	f_{j}[n]=\begin{cases}
	1, &\text{for} \ n=175,525,875,...,4025 \\
	0, &\text{otherwise}
	\end{cases}
	\ \ \ \ \ \forall j\in\mathcal{J}_1 \\
	f_{j}[n]=\begin{cases}
	1, &\text{for} \ n=0,175,350,525,700,875,...,4025 \\
	0, &\text{otherwise}
	\end{cases}
	\ \ \ \ \ \forall j\in\mathcal{J}_2
	\end{align*}
\end{subequations}

A single-subconstellation system is assumed with the corresponding seed satellite orbital elements vector: $\textbf{\oe}_0=[83/6,0,99.2\degree,0\degree,0\degree,0\degree]^T$. This orbit corresponds to an altitude of $\SI{946.7}{km}$. The period of repetition of this orbit is $T_{\text{r}}=\SI{5.184e05}{s}$, which is six days. The system must satisfy:

\begin{equation*}
\renewcommand{\arraystretch}{1.2}
\begin{bmatrix}
\bm{V}_{0,1} \\
\bm{V}_{0,2} \\
\vdots \\
\bm{V}_{0,86}
\end{bmatrix}
\bm{x}
\ge
\begin{bmatrix}
\bm{f}_1 \\
\bm{f}_2 \\
\vdots \\
\bm{f}_{86} \\
\end{bmatrix}
\end{equation*}

{\parindent0pt
	where the dimension of this inequality is $(361200\times4200)\cdot (4200\times1)\ge(361200\times1)$.
}

The results show that the quasi-symmetric constellation is composed of 96 satellites, whereas the BILP constellation is composed of 29 satellites. Comparing the computational cost, the quasi-symmetric method took \SI{486.7}{s}, whereas the BILP method took only \SI{7.7}{s}. This shows a significant improvement of the BILP method in terms of both the number of satellites and the computational time with respect to the quasi-symmetric constellation. The quasi-symmetric method is performing poorly because we need a large number of satellites if the symmetric pattern is used. This factor, together with the large numbers of target points and time steps, makes the iterative process in the quasi-symmetric method inefficient. The BILP method, instead, identifies the asymmetric optimal solution with a significantly smaller number of satellites. The low computational cost for the BILP method is due to the BILP solver, Gurobi in our case; the problem structure allows Gurobi to perform an efficient presolve procedure, resulting in a short optimization time.

\begin{figure}[H]
	\centering
	\renewcommand{\arraystretch}{0.8}
	\begin{tabular}{cccc}
		\includegraphics[height=4cm]{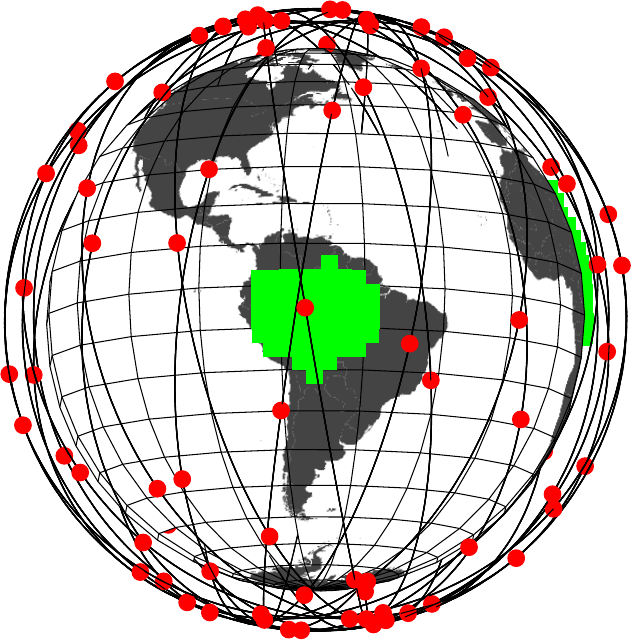} &
		\includegraphics[height=4cm]{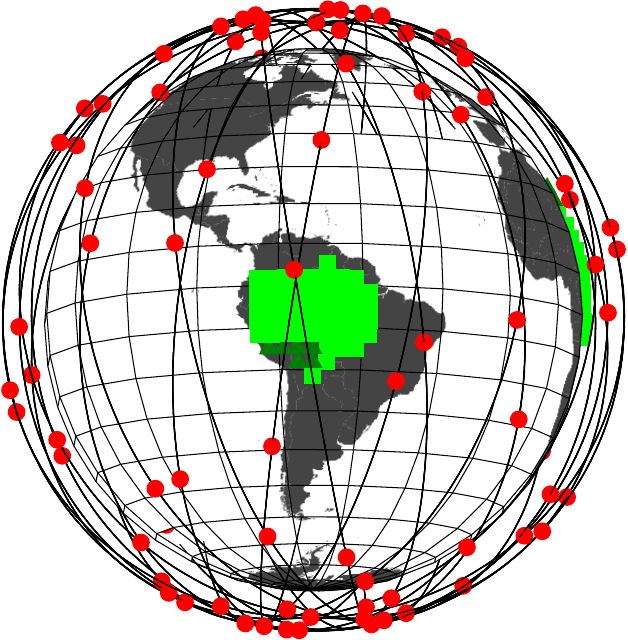} &
		\includegraphics[height=4cm]{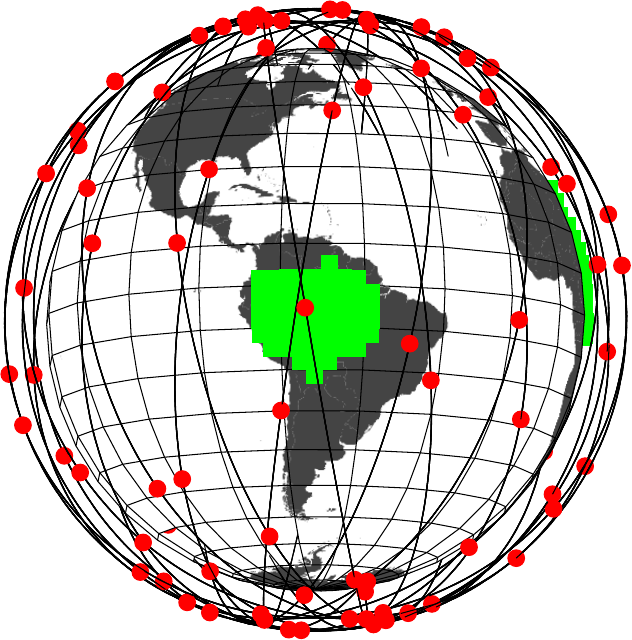} \\
		\textbf{(a)} Quasi-symmetric: $n=0$ & \textbf{(b)} Quasi-symmetric: $n=88$ & \textbf{(c)} Quasi-symmetric: $n=175$ \\
		req:$f[0]=0$; result cov: \SI{100}{\%} & req:$f[88]=0$; result cov: \SI{89}{\%} & req:$f[175]=1$; result cov: \SI{100}{\%}\\ [8pt] 
	\end{tabular}
	\begin{tabular}{cccc}
		\includegraphics[height=4cm]{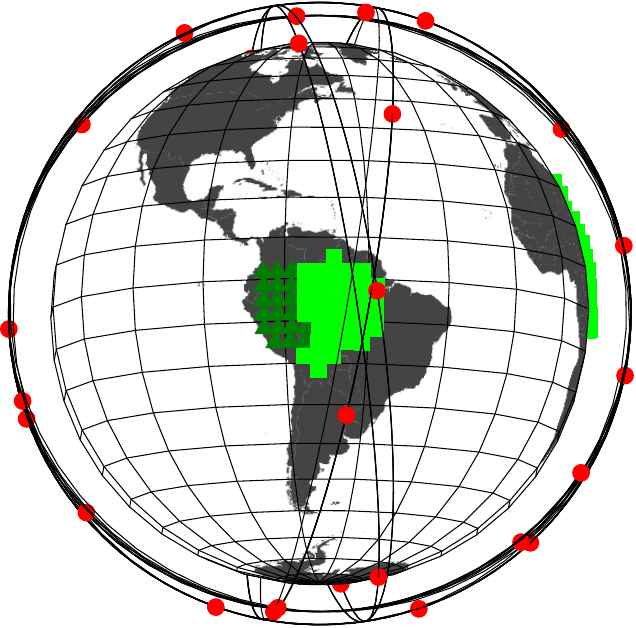} &
		\includegraphics[height=4cm]{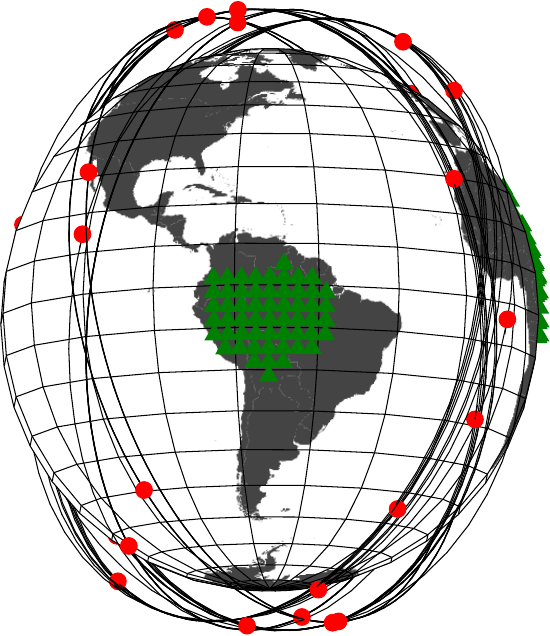} &
		\includegraphics[height=4cm]{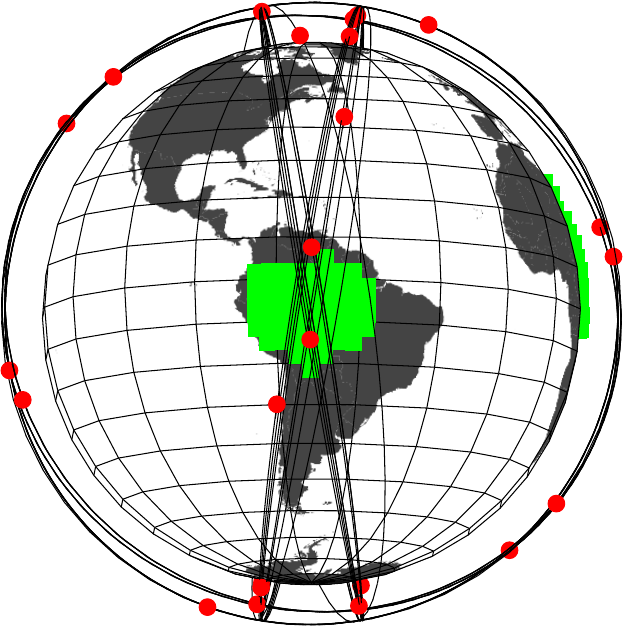} \\
		\textbf{(d)} BILP: $n=0$ &\textbf{(e)} BILP: $n=88$ & \textbf{(f)} BILP: $n=175$ \\ req:$f[0]=0$; result cov: \SI{66}{\%}  & req:$f[88]=0$; result cov: \SI{0}{\%} & req:$f[175]=1$; result cov: \SI{100}{\%} \\[8pt]
	\end{tabular}
	\caption{Example 4: Coverage over the Amazon river basin; select snapshots are shown at $n=0,88,175$ (ECI frame). (a), (b), and (c) are the snapshots of the quasi-symmetric constellation and (d), (e), and (f) are the snapshots of the BILP constellation. At each $n$, targets that have satellite visibility are shown in light green squares and targets that do not have satellite visibility are shown in dark green triangles. ``req'' indicates the coverage requirement, and ``result cov'' is the actual coverage performance of the solution. For example, when the requirement $f[n]=1$, the coverage has to be 100\% (i.e., at least one satellite is visible from all target points in the area). It can be seen that the BILP method takes advantage of asymmetry and satisfies the coverage requirements with fewer satellites.}
	\label{fig:example4_3d_amazon}
\end{figure}

\begin{figure}[H]
	\centering
	\renewcommand{\arraystretch}{0.8} 
	\begin{tabular}{cccc}
		\includegraphics[height=4cm]{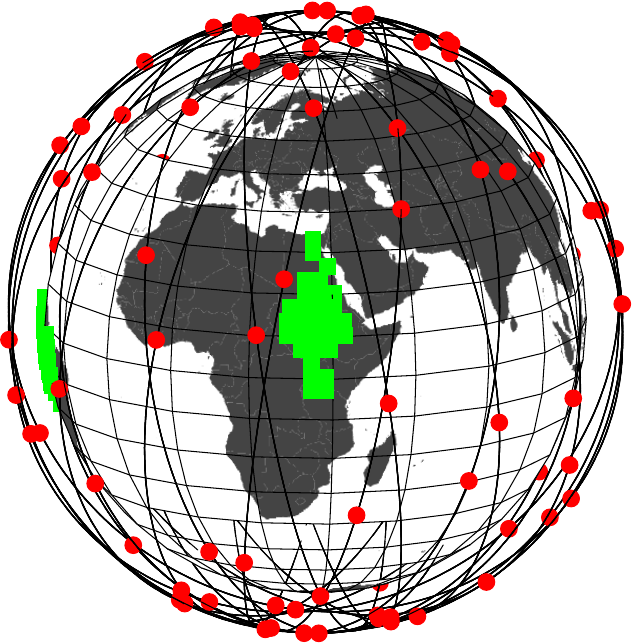} &
		\includegraphics[height=4cm]{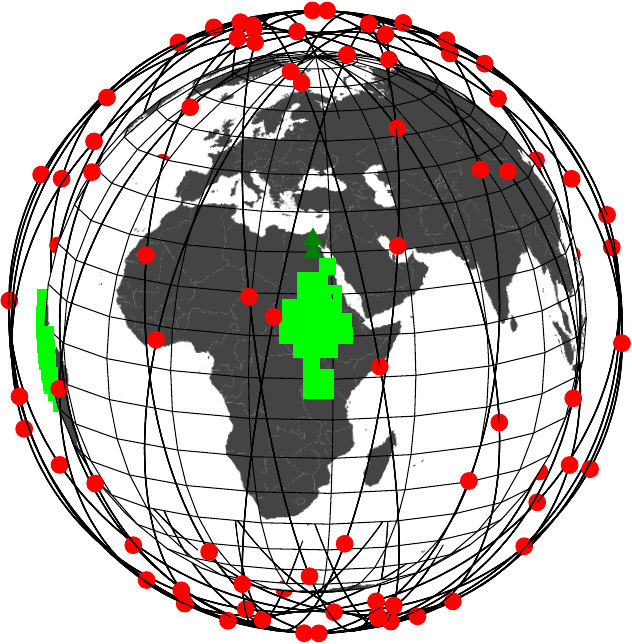} &
		\includegraphics[height=4cm]{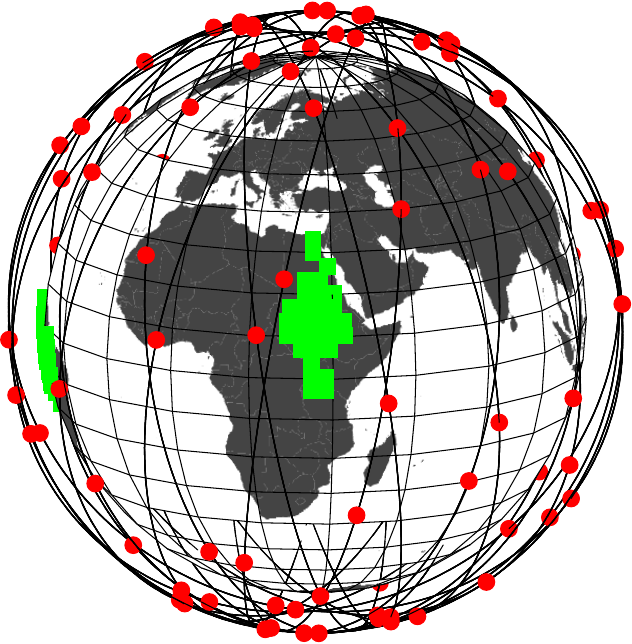} \\
		\textbf{(a)} Quasi-symmetric: $n=0$ & \textbf{(b)} Quasi-symmetric: $n=88$ & \textbf{(c)} Quasi-symmetric: $n=175$ \\
		req:$f[0]=1$; result cov: \SI{100}{\%} & req:$f[88]=0$; result cov: \SI{93}{\%} & req:$f[175]=1$; result cov: \SI{100}{\%}  \\[8pt]
	\end{tabular}
	\begin{tabular}{cccc}
		\includegraphics[height=4cm]{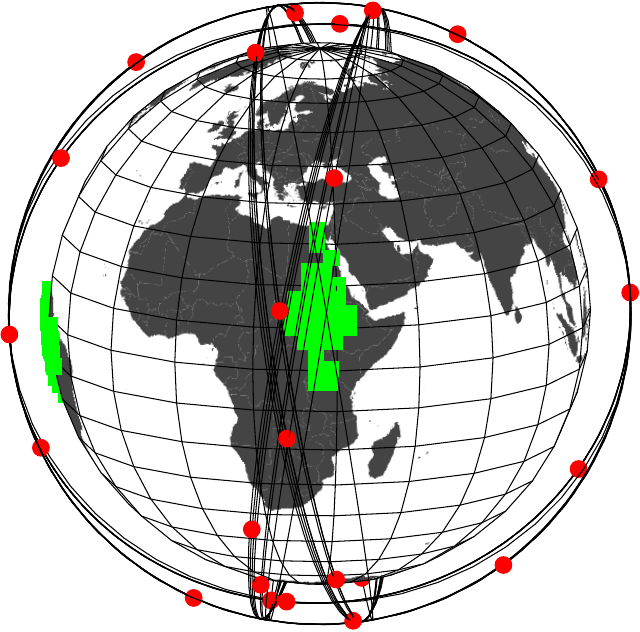} &
		\includegraphics[height=4cm]{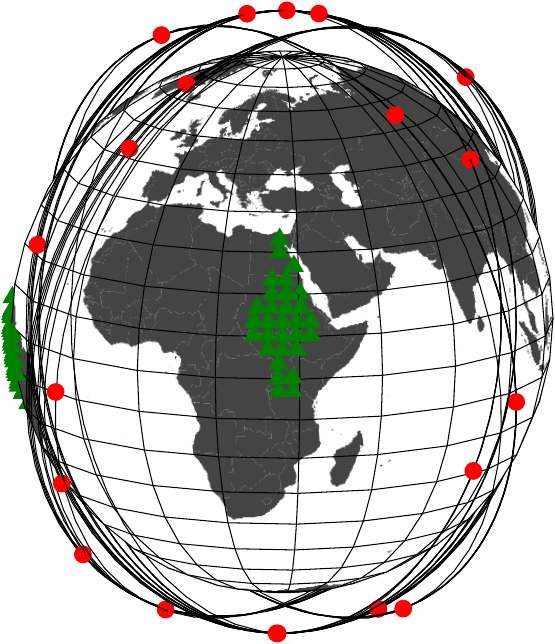} &
		\includegraphics[height=4cm]{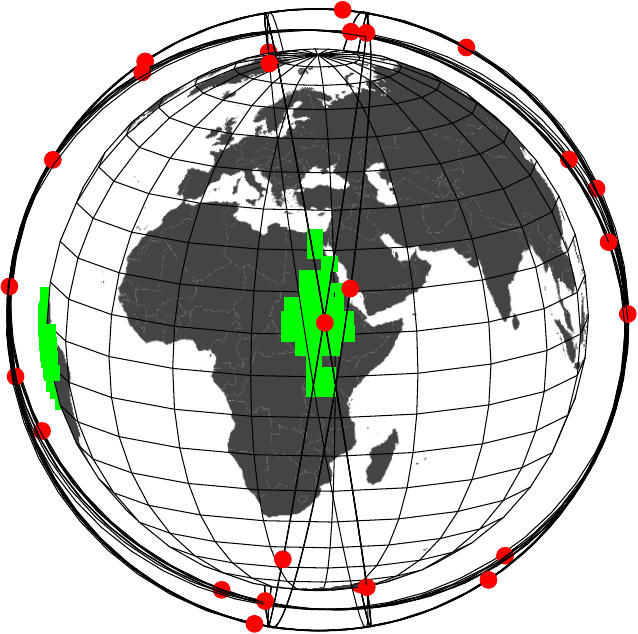} \\
		\textbf{(d)} BILP: $n=0$ &\textbf{(e)} BILP: $n=88$ & \textbf{(f)} BILP: $n=175$ \\
		req:$f[0]=1$; result cov: \SI{100}{\%} & req:$f[88]=0$; result cov: \SI{0}{\%} & req:$f[175]=1$; result cov: \SI{100}{\%} \\[8pt]
	\end{tabular}
	\caption{Example 4: Coverage over the Nile river basin; select snapshots are shown at $n=0,88,175$ (ECI frame). (a), (b), and (c) are the snapshots of the quasi-symmetric constellation and (d), (e), and (f) are the snapshots of the BILP constellation. At each $n$, targets that have satellite visibility are shown in light green squares and targets that do not have satellite visibility are shown in dark green triangles. ``req'' indicates the coverage requirement, and ``result cov'' is the actual coverage performance of the solution. For example, when the requirement $f[n]=1$, the coverage has to be 100\% (i.e., at least one satellite is visible from all target points in the area). It can be seen that the BILP method takes advantage of asymmetry and satisfies the coverage requirements with fewer satellites.}
	\label{fig:example4_3d_nile}
\end{figure}

Fig.~\ref{fig:example4_3d_amazon} and Fig.~\ref{fig:example4_3d_nile} show the select snapshots of both the quasi-symmetric constellation and the BILP constellation in chronological order over the Amazon and Nile river basins, respectively. (Due to the large number of satellites, the resulting constellation pattern vector is omitted.) As expected, both constellations provide simultaneous access to the target points when needed ($n=175$ for Amazon river basin; $n=0,175$ for the Nile river basin), satisfying the strict revisit time requirements. It can be seen that, while the quasi-symmetric method satisfies the coverage requirements with a (quasi-)symmetric constellation pattern, the BILP method takes advantage of the asymmetry and satisfies the same requirements with fewer satellites.

\subsection{Example 5. A System of Multiple Sub-Constellations over Multiple Target Points}
We consider a most general case that only the BILP method can solve: a system of multiple sub-constellations over multiple target points. In this example, two target points are considered in the target point set $\mathcal{J}=\{(\phi=64.14\degree \text{N},\lambda=21.94\degree \text{W}), (\phi=19.07\degree \text{N},\lambda=72.87\degree \text{E})\}$: Reykjav\'ik, Iceland ($j=1$) and Mumbai, India ($j=2$). The minimum elevation angle for each target point is: $\varepsilon_{1,\text{min}}=15\degree$ and $\varepsilon_{2,\text{min}}=10\degree$. The objective is to achieve single-fold continuous coverage over all target points ($\bm{f}_{j}=\bm{1}, \ \forall j\in{\mathcal{J}}$).

Two sub-constellations are considered: $\textbf{\oe}_0^{(1)}=[8/1,0,70\degree,0\degree,0\degree,0\degree]^T$ (an altitude of $\SI{4149.2}{km}$) and $\textbf{\oe}_0^{(2)}=[6/1,0,47.915\degree,0\degree,0\degree,0\degree]^T$ (an altitude of $\SI{6380.3}{km}$). The length of vectors is selected, $L=717$, such that the time step is approximately $t_{\text{step}}\approx\SI{120}{s}$. The period of repetitions for these sub-constellations are identical, $T_{\text{r}}^{(1)}=T_{\text{r}}^{(2)}\approx\SI{86024}{s}$, hence making two sub-constellations synchronous. Note that, even though we are using two sub-constellations for disconnected regions of interest, the sub-constellations are not defined one per region of interest; instead, they are used together to satisfy both demands in an optimal way. The goal of the BILP method is to optimize $\bm{x}^{(1)}$ and $\bm{x}^{(2)}$ concurrently such that the system satisfies the augmented linear condition:

\begin{equation*}
\renewcommand{\arraystretch}{1.2}
\begin{bmatrix}
\bm{V}_{0,1}^{(1)} & \bm{V}_{0,1}^{(2)} \\
\bm{V}_{0,2}^{(1)} & \bm{V}_{0,2}^{(2)} \\
\end{bmatrix}
\begin{bmatrix}
\bm{x}^{(1)} \\
\bm{x}^{(2)} \\
\end{bmatrix}
\geq
\begin{bmatrix}
\bm{f}_1 \\
\bm{f}_2 \\
\end{bmatrix}
\Leftrightarrow
\{\bm{V}_{0,1}^{(1)}\bm{x}^{(1)}+\bm{V}_{0,1}^{(2)}\bm{x}^{(2)}\geq\bm{f}_1, \bm{V}_{0,2}^{(1)}\bm{x}^{(1)}+\bm{V}_{0,2}^{(2)}\bm{x}^{(2)}\geq\bm{f}_2\}
\end{equation*}

The following optimal constellation pattern vectors are obtained:
\begin{subequations}
	\begin{align*}
	x^{(1)\ast}[n]&=\begin{cases}
	1, &\text{for} \ n=65,144,285,361 \\
	0, &\text{otherwise}
	\end{cases} \\
	x^{(2)\ast}[n]&=\begin{cases}
	1, &\text{for} \ n=208,428,523,608,634,702 \\
	0, &\text{otherwise}
	\end{cases}
	\end{align*}
\end{subequations}

The number of satellites is 4 for the first sub-constellation and 6 for the second; 10 in total. The computational time was \SI{5298.7}{s}. 

\begin{figure}[H]
	\centering
	\begin{subfigure}[h]{0.495\linewidth}
		\centering
		\includegraphics[trim=43 10 45 15,clip,width=\linewidth]{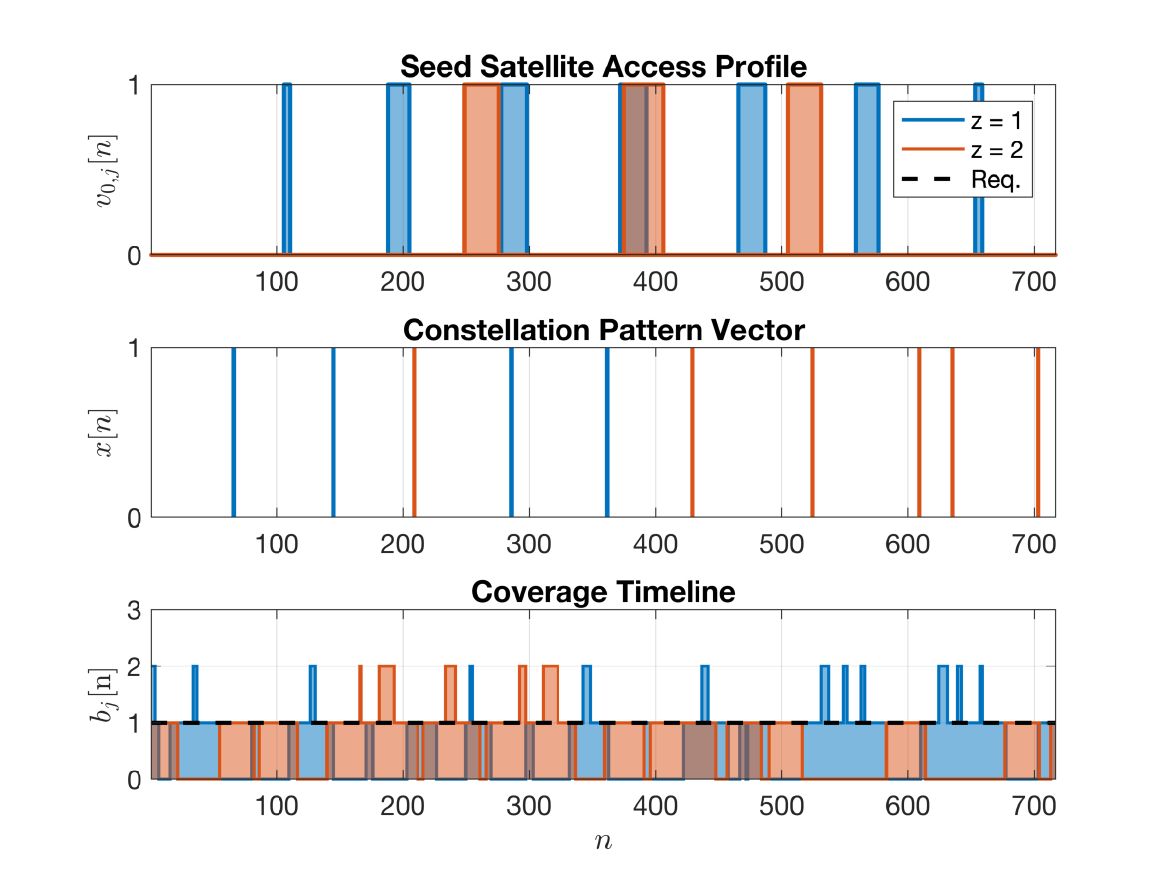}
		\caption{Individual contribution over Reykjav\'ik, Iceland}
	\end{subfigure}
	\begin{subfigure}[h]{0.495\linewidth}
		\centering
		\includegraphics[trim=43 10 45 15,clip,width=\linewidth]{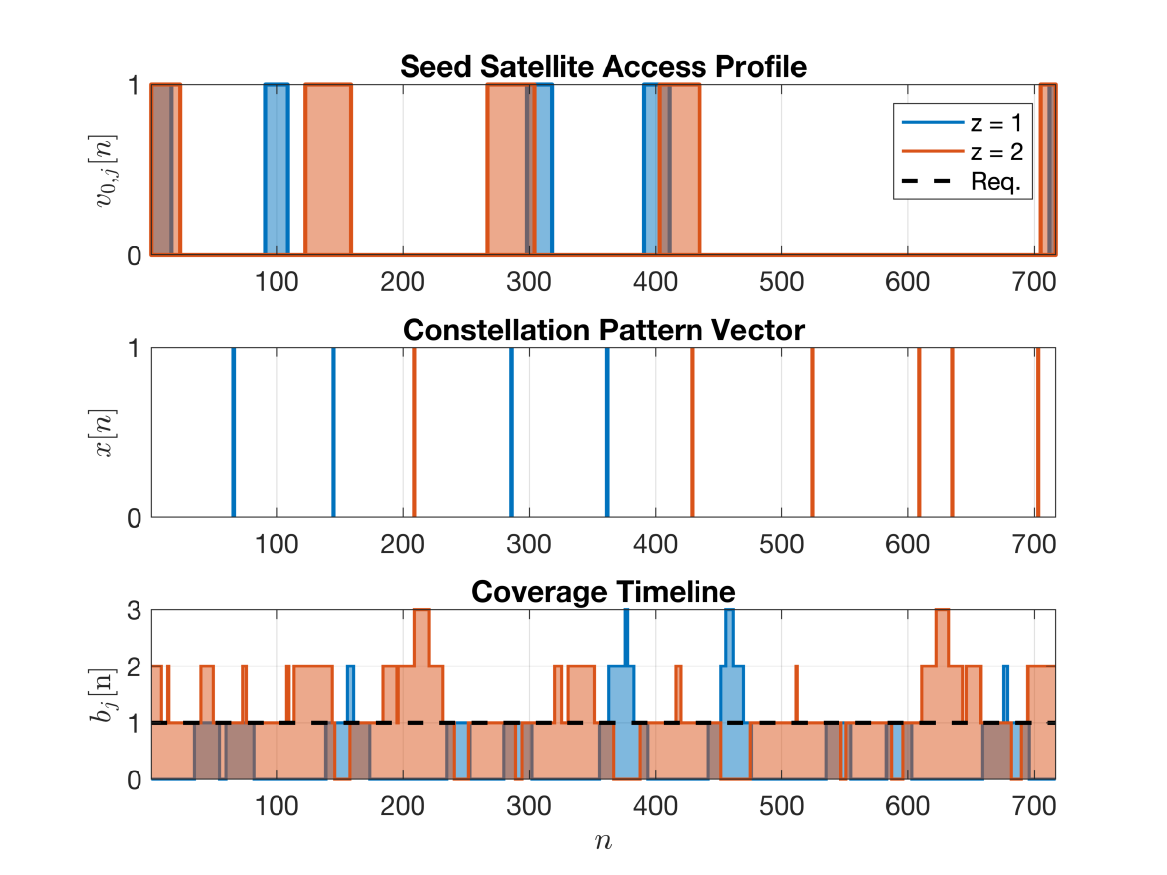}
		\caption{Individual contribution over Mumbai, India}
	\end{subfigure}
	\caption{Example~5: The APC decomposition}
	\label{fig:example5_apc}
\end{figure}

Fig.~\ref{fig:example5_apc} illustrates the benefit of the BILP method. Individually, $z=1$ sub-constellation provides \SI{53.7}{\%} and \SI{37.1}{\%} coverage over $j=1$ and $j=2$, respectively and $z=2$ sub-constellation provides \SI{65.0}{\%} and \SI{87.0}{\%} coverage over $j=1$ and $j=2$, respectively. No individual sub-constellation alone provides complete continuous coverage over any target point. The BILP method concurrently optimizes $\bm{x}^{(1)}$ and $\bm{x}^{(2)}$ such that the continuous coverage over the whole target set $\mathcal{J}$ is achieved while minimizing the total number of satellites from two sub-constellations. Note that the constellation pattern vectors, $\bm{x}^{(1)\ast}$ and $\bm{x}^{(2)\ast}$, are identical in both sub-figures of Fig.~\ref{fig:example5_apc}.

\begin{figure}[H]
	\centering
	\begin{subfigure}[h]{0.4\linewidth}
		\centering
		\includegraphics[height=5cm]{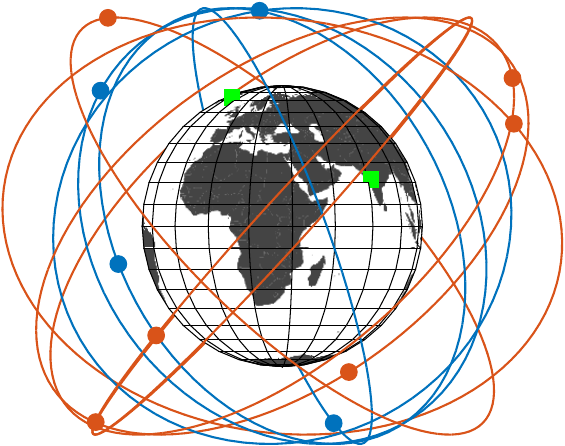}
		\caption{View from side}
	\end{subfigure}
	\begin{subfigure}[h]{0.4\linewidth}
		\centering
		\includegraphics[height=5cm]{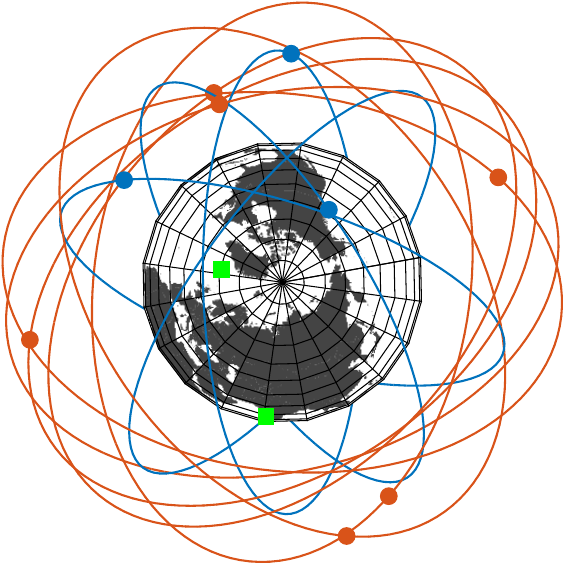}
		\caption{View from the North Pole}
	\end{subfigure}
	\caption{Example 5: 3D view of generated constellation at $n=0$ (ECI frame)}
	\label{fig:example5_3d}
\end{figure}

The optimized two-subconstellation system is shown in Fig.~\ref{fig:example5_3d}. The sub-constellation ($z=1$) colored in blue (lower altitude) is composed of four satellites while the sub-constellation ($z=2$) colored in red (higher altitude) is composed of six satellites for a total of ten satellites.

Finally, to show the effectiveness of having the sub-constellations, corner cases are evaluated considering each individual sub-constellation separately. The results indicate that, under the same setting, using only the sub-constellation 1 results in 11 satellites, and using only sub-constellation 2 also results in 11 satellites. This particular case demonstrates that through the use of multiple sub-constellations, one can reduce the minimum satellites required from 11 to 10 by enlarging the design space. Also, it is worth mentioning that the BILP method can still lead to an optimal solution even for the cases where only part of the sub-constellation sets is used in the optimal pattern.

\section{Conclusion} \label{sec:conclusion}
A semi-analytical approach to optimally design a regional coverage satellite constellation pattern is proposed. By treating the seed satellite access profile and the constellation pattern vector as discrete-time signals, a circular convolution between them creates the coverage timeline. We refer to this formulation as the APC decomposition of the satellite constellation system. This formulation is used to derive a set of satellite constellation pattern design methods that take a seed satellite access profile and a coverage requirement as their inputs and output the minimum number of satellites required to satisfy the coverage requirement. Two satellite constellation pattern design methods are introduced: the baseline quasi-symmetric method and the more general BILP method. The baseline quasi-symmetric method enforces the conventional assumption of symmetry in the constellation pattern and solves for the minimum number of satellites required in the system by incrementally increasing $N$ until the coverage requirement is satisfied. In contrast, the new and more general BILP method solves for constellation pattern vector $\bm{x}$ where $N$ and their temporal locations can be deduced by solving a binary integer linear programming problem. Our analysis shows that, while the quasi-symmetric method can be efficient when we can satisfy the coverage requirements with a small number of satellites in a symmetric pattern (e.g., continuous polar coverage), the BILP method always outputs optimal satellite constellation patterns that the baseline method may miss. Furthermore, the BILP method is applicable to the problems that the quasi-symmetric method cannot solve (e.g., the case with multiple sub-constellations).

Our ideas respond to the several design features that can reinforce the utility of regional constellations: multiple target points, complex time-varying and spatially-varying requirements, and multiple sub-constellations. The developed circular convolution formulation allows linearity in both the multiple target points direction and multiple sub-constellations direction via matrix augmentation. A user can design (1) a single constellation system that simultaneously satisfies the complex coverage requirement of area targets composed of multiple target points, (2) a system of multiple sub-constellations that satisfies the complex coverage requirement of a single target point, or (3) a combination of both. These design features are demonstrated via a series of illustrative examples in Section~\ref{sec:illustrative_examples}. The resulting general constellation pattern design approach can be integrated with existing orbital characteristics design methods and launch/mission constraints to help future satellite constellation designers rigorously achieve optimal constellation designs.

Despite the demonstrated effectiveness of the proposed approach, there are some possible directions for future work to improve it further. The first potential direction is related to the computational time. Due to the nature of the discretization, obtaining a high-fidelity solution computed with fine time discretization would require a large-sized problem and thus a long computational time. To make the method computationally more scalable, approximation algorithms or heuristics methods can be developed to retrieve feasible, yet potentially suboptimal, solutions in a relatively short amount of time. Furthermore, the proposed method only considers the $J_2$ effect as the disturbance and assumes that the spacecraft has the maneuvering capability to cancel out other disturbances. This assumption is reasonable for the proposed method to be used for a high-level constellation pattern design purpose, but it can be improved for higher-fidelity modeling. Finally, this constellation pattern design method requires the seed satellite orbital elements as its input. While Appendix~C shows one example process of integrating the proposed approach into the constellation design practice, further investigation can be performed to ensure an efficient and effective integration. 

\section*{Appendix~A. Expanded Ground Track View} \label{app:egtv}
\begin{figure}[htb]
	\centering\includegraphics[width=\textwidth]{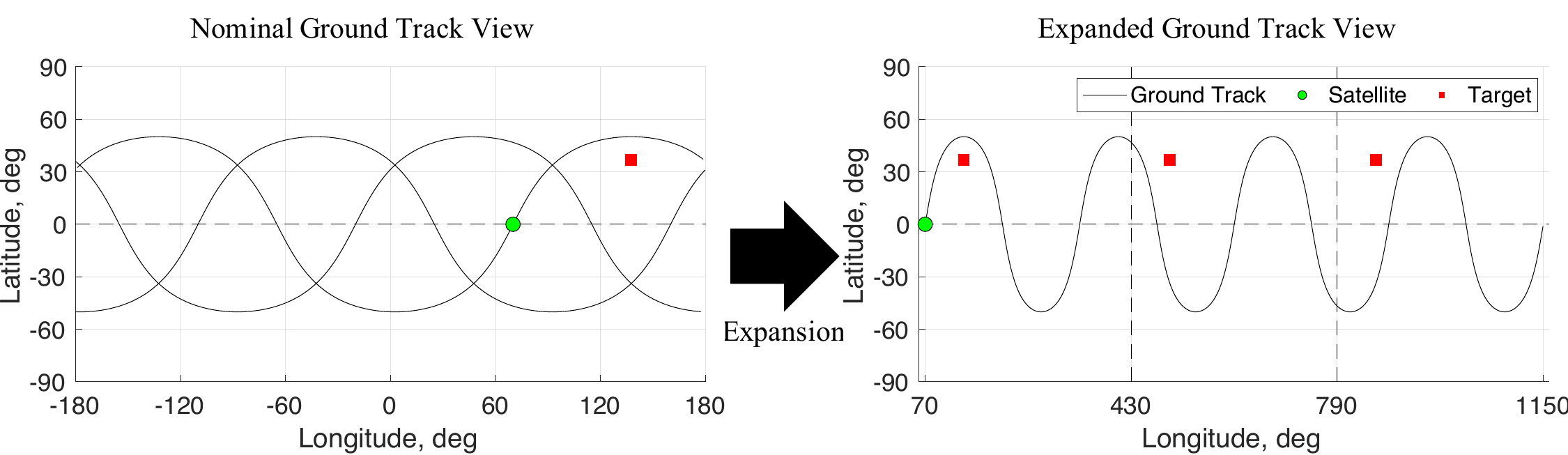}
	\caption{Full expansion of a ground track of $\textbf{\oe}_0=[4/1,0,50\degree,0\degree,350.2\degree,0\degree]$ (J2000)}
	\label{fig:egtv}
\end{figure}

The \textit{expanded ground track view} spatially expands an ordinary ground track of a satellite and visualizes its ground track relative to the area of interest throughout the simulation period $T_{\text{sim}}$. The area of interest and its mirrored images are positioned throughout the plot (the red squares in Fig.~\ref{fig:egtv}) to provide spatial references. The expanded ground track view is especially useful when visualizing and correlating the access profile and the actual satellite ground track.

The following properties of the expanded ground track view are formalized for the repeating ground track with the period ratio of $\tau=N_{\text{P}}/N_{\text{D}}$.

\begin{enumerate}
    \item The magnitude of the longitudinal angular displacement of the expanded ground track is $360\lvert N_{\text{P}}-N_{\text{D}}\rvert$ degrees for prograde orbits or $360(N_{\text{P}}+N_{\text{D}})$ degrees for retrograde orbits \cite{kim1997theory}. Here, the longitudinal angular displacement of the expanded ground track is defined as the total angular displacement required to repeat the ground track, measured along the axis of longitude in the direction of the satellite's motion.
	\item The mirrored images of the area of interest are separated by 360 degrees.
\end{enumerate}

\section*{Appendix~B. Derivation of the Coverage Timeline} \label{app:derivation}
To prove the circular convolution phenomenon, show that Eq.~\eqref{eq:b} is identical to Eq.~\eqref{eq:circulantX}. Begin by expanding Eq.~\eqref{eq:b}, which is the summation of all access profiles:

\begin{equation}
\label{eq:ba}
b_j[n]=v_{1,j}[n]+v_{2,j}[n]+\cdots+v_{N,j}[n]
\end{equation}

Each term of Eq.~\eqref{eq:ba} can be represented as a multiple of $v_{0,j}[n]$ and permutation matrix $\bm{P}_\pi^{n_k}$ due to the cyclic property of the assumed formulation. Recalling the definition from Eq.~\eqref{eq:bb}:

\begin{equation*}
v_{k,j}[n]=\bm{P}_{\pi}^{n_k}v_{0,j}[n]
\end{equation*}

{\parindent0pt
	where $\bm{P}_\pi$ is a permutation matrix with the dimension $(L \times L)$ shown below. Note that $\bm{I}=\bm{P}_\pi^0=\bm{P}_\pi^L$.
}

\begin{equation}
\bm{P}_\pi=
\renewcommand{\arraystretch}{0.8}
\begin{bmatrix}
0 & 0 & 0 & \cdots &1 \\
1 & 0 & 0 & \cdots & 0 \\
0 & 1 & 0 & \ddots & \vdots \\
\vdots & \vdots & \ddots & \ddots & 0 \\
0 & 0 & \cdots & 1 & 0
\end{bmatrix}
\end{equation}

Substituting Eq.~\eqref{eq:bb} into Eq.~\eqref{eq:ba}, we get the following equation:

\begin{equation}
\label{eq:bc}
b_j[n]=(\bm{P}_\pi^{n_1}+\bm{P}_\pi^{n_2}+\cdots+\bm{P}_\pi^{n_N})v_{0,j}[n]
\end{equation}

Eq.~\eqref{eq:bc} is a superposition of cyclically shifted access profiles referenced to a seed satellite access profile. Here, $n_k$ denotes the index of the relative time shift of the $k$th access profile with respect to the seed satellite access profile. Instead of only indicating the indices where only access profiles exist, one can generalize this to all time steps $n\in\{0,...,L-1\}$ following the definition of the constellation pattern vector in Eq.~\eqref{eq:x}. Hence, Eq.~\eqref{eq:bc} can be further deduced as:

\begin{equation}
\label{eq:bd}
b_j[n]=\bigg(x[0]\bm{P}_\pi^{0}+x[1]\bm{P}_\pi^{1}+\cdots+x[L-1]\bm{P}_\pi^{L-1}\bigg)v_{0,j}[n]
\end{equation}

The terms within parentheses in Eq.~\eqref{eq:bd} is identical to the alternative analytical definition of the circulant matrix:

\begin{equation}
\label{eq:be}
\bm{X} \triangleq x[0]\bm{I}+x[1]\bm{P}_\pi^1+\cdots+x[L-1]\bm{P}_\pi^{L-1}
\end{equation}

Finally, substituting Eq.~\eqref{eq:be} into Eq.~\eqref{eq:bd}, we get:

\begin{equation}
\label{eq:bf}
\bm{b}_j=\bm{X}\bm{v}_{0,j}
\end{equation}

Using the commutative property of the circular convolution operator, Eq.~\eqref{eq:bf}:

\begin{equation}
\bm{b}_j=\bm{V}_{0,j}\bm{x}
\end{equation}
where 
\begin{equation*}
V_{0,j}[\alpha,\beta] = v_{0,j}[(\alpha-\beta) \bmod L]
\end{equation*}
as defined in Eq.~\eqref{eq:vdef}.

This is identical to the definition of the circular convolution in Eq.~\eqref{eq:circulantX}, thereby proving the circular convolutional nature of the formulation under the aforementioned assumptions.

\section*{Appendix~C. Integrating the Developed Method into Constellation Design Process}

This appendix introduces an example approach to integrate the developed method into the satellite constellation design process. As discussed earlier, the developed satellite constellation pattern design method needs the seed satellite orbital elements as its input. In this appendix, we introduce an approach to efficiently integrate the determination of the seed satellite orbital elements $\textbf{\oe}_0$ and the determination of the constellation pattern $\bm{x}$ (i.e., the developed method).

First, note that although $\textbf{\oe}_0$ contains six orbital elements ($\tau,e,i,\omega,\Omega_0,M_0$), we only have five degrees of freedom. The initial mean anomaly of the seed satellite $M_0$ can be set to zero without loss of generality. This is because, as shown in Eq.~\eqref{eq:set}, $\Omega_0$ and $n_k$ can be chosen such that any solution with an arbitrarily chosen $M_0$ can be converted into an equivalent solution with $M_0=0\degree$\footnote{Strictly speaking, there are only a finite number of possible discrete values for $M_0$ due to the discretization used in this problem.}.

The design space of the remaining five orbital elements can be narrowed down even further by considering the launch and mission requirements. As an example, we consider the case used in Example 2 in Section~\ref{sec:illustrative_examples} and provide a walk-through process.
\begin{enumerate}
    \item Suppose there is demand for increased communications capacity (i.e., increased satellite diversity) during a particular time interval of a day that repeats daily (e.g., Internet rush hour) over Atlanta, Georgia ($\{(\phi=34.75\degree \text{N},\lambda=84.39\degree \text{W})\}$). Translating this demand, the time-varying coverage requirement $\bm{f}$ is derived (see Example 2). The communications quality-of-service requirement further enforces consistency in data round-trip latency throughout the mission duration; hence, a circular orbit is desired. The period ratio and the minimum elevation angle are assumed to be derived a priori based on mission-related requirements: $\tau=12/1$ and $\varepsilon_{\text{min}}=5\degree$.
    \item Based on the set of mission requirements and parameters ($T_\text{r}=\SI{86400}{s}$, $e=0$, and $\tau=12/1$), the inclination of the orbit is readily derived, which is approximately 102.9\degree. Note that since the repeat period $T_\text{r}$ is exactly given together with $\tau$ and $e$, there is no degree of freedom for trading off the altitude and the inclination. In this case, since the repeat period is exactly $\SI{86400}{s}$, the orbit needs to be a repeating sun-synchronous orbit.
    \item At this point, the only leftover variable is $\Omega_0$, which dictates the shift of the common ground track along the longitudinal direction. The RAAN of the seed satellite $\Omega_0$ can be determined either by an analytical heuristics method or by a numerical optimization.
    \begin{enumerate}[label=\alph*)]
        \item An analytical heuristic approach can determine $\Omega_0$ such that the common ground track is symmetric about the longitude of the target point (see Fig.~\ref{fig:symmetry}). Solving for the corresponding RAAN value yields $\Omega_0=98.3\degree$. Note that another symmetry exists further offsetting $\Omega_0$ value.
        \item A single-variable optimization can be performed to determine the value of $\Omega_0$. Ideally, we prefer to use the number of satellites as the metric, but this cannot be evaluated without $\bm{x}$. Instead, an effective metric can be the coverage over the area of interest. Note that the values of $\textbf{\oe}_0$ maximizing the coverage does not necessarily lead to a minimum number of satellites, but as shown later, it is a good approximation to use.
    \end{enumerate}
    
    \begin{figure}[htb]
	\centering\includegraphics[width=0.55\textwidth]{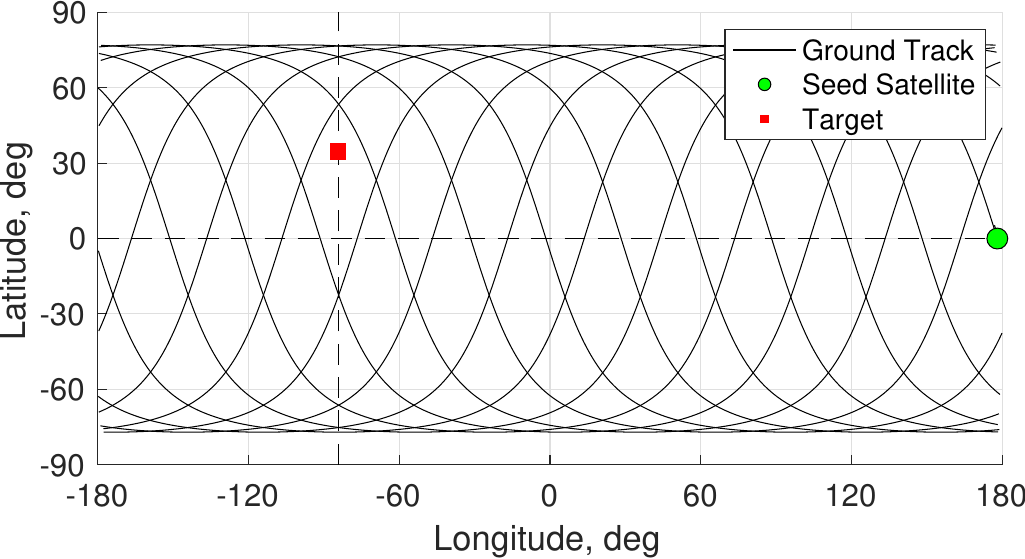}
	\caption{Alignment of the ground track such that it is symmetric about the longitude of the target}
	\label{fig:symmetry}
	\end{figure}
    
    \item Using the obtained seed satellite orbital elements, the optimization of the constellation pattern vector can be performed following the APC-based methods developed in this paper.
\end{enumerate}

As we evaluate the efficiency of the developed integrated heuristics and BILP methods, we compare them against a more straightforward approach, where both $\textbf{\oe}_0$ and $\bm{x}$ are optimized as variables simultaneously against the objective function of the number of satellites. In fact, this formulation is the most direct representation of our goal; however, since it is a mixed-integer nonlinear optimization problem, we cannot leverage the developed method in this paper and therefore can only use generic inefficient solvers (e.g., genetic algorithm). Here, we aim to show that, by incorporating the developed method into this process, we can achieve a much better performance than this classical integrated method.

In Table~\ref{tab:comparison}, Method 1 refers to the heuristics approach that finds $\Omega_0$ using symmetry, which is the actual method used in Example~2. Method 2 refers to the two-stage optimization where the first stage is the metaheuristics optimization of $\Omega_0$ and the second stage is the BILP optimization of $\bm{x}$. Lastly, Method 3 is the simultaneous optimization of both $\Omega_0$ and $\bm{x}$ via metaheuristics optimization. For Methods 2 and 3, a genetic algorithm (GA) by MATLAB is used with the default settings.

\begin{table}[H]
	\fontsize{10}{10}\selectfont
	\caption{Comparison of different methods for integrated optimization}
	\centering
	\begin{tabular}{c c c}
		\hline
		\hline
		Method & Number of Satellites & Computational Time, s \\
		\hline
		1 & 24 & 3712.0 \\
		\multirow{2}{*}{2} & \multirow{2}{*}{25} & Stage 1: 477.1 \\
		& & Stage 2: 2086.0 \\
		\multirow{4}{*}{3} & 75 & 2482.1 (Population: 100) \\
		& 36 & 5854.7 (Population: 200) \\
		& 32 & 10635.1 (Population: 300) \\
		\hline
		\hline
	\end{tabular}
	\label{tab:comparison}
\end{table}

The results show that both Methods 1 and 2 are effective in finding the optimal solution; the only difference in these two methods is in the optimization of $\Omega_0$. On the other hand, Method 3 requires longer computational time, while only showing poor results. These results demonstrate the utility of the developed method when integrated into the satellite constellation design process.

\section*{Appendix~D. Derivation of the RAAN phasing}
This appendix derives \[\Omega_k=n_k\frac{2\pi N_\text{D}}{L}+\Omega_0\] in Eq.~\eqref{eq:setb}. Define $\Delta\Omega=(\Omega_k-\Omega_0)/n_k$. Our goal is to prove $\Delta\Omega=2\pi N_\text{D}/L$.

This expression comes from Fig.~\ref{fig:pattern}. In order to achieve a constellation that separates away from each other by $t_\text{step}$ over a common ground track, $\Delta\Omega$ needs to be defined as the difference between Earth's rotation and the angular displacement due to the RAAN precession during a time interval $[0,t_{\text{step}}]$. More specifically,

\begin{equation}
\label{eq:delta_omega}
{\Delta\Omega}=(\omega_{\oplus}-\dot{\Omega}){t_{\text{step}}}
\end{equation}
Since $t_{\text{step}}=T_{\text{r}}/L$, substituting in Eq.~\eqref{eq:period_of_repetition} yields $t_{\text{step}}=N_\text{D}T_\text{G}/L$. Plugging this into Eq.~\eqref{eq:delta_omega}, we get:

\begin{equation}
\Delta\Omega = (\omega_{\oplus}-\dot{\Omega})\frac{N_\text{D}T_\text{G}}{L}
\end{equation}
Since, $T_\text{G}=2\pi/(\omega_{\oplus}-\dot{\Omega})$ (Eq.~\eqref{eq:tg}), we get:

\begin{equation}
\Delta\Omega = \frac{2\pi N_\text{D}}{L}
\end{equation}

\section*{Acknowledgments}
This research is supported by the Advanced Technology R\&D Center at Mitsubishi Electric Corporation. The first author would like to acknowledge additional support from the National Science Foundation. This material is based upon work supported by the National Science Foundation Graduate Research Fellowship Program under Grant No. DGE–1650044. Any opinions, findings, and conclusions or recommendations expressed in this material are those of the author(s) and do not necessarily reflect the views of the National Science Foundation. The authors would like to thank Onalli Gunasekara, Hao Chen, and Dr. Robert Griffin for their editorial review and thoughtful suggestions for improvement.

\pagebreak

\bibliography{references}

\begin{thebibliography}{43}
\newcommand{\enquote}[1]{``#1''}
\providecommand{\natexlab}[1]{#1}
\providecommand{\url}[1]{\texttt{#1}}
\providecommand{\urlprefix}{URL }
\expandafter\ifx\csname urlstyle\endcsname\relax
  \providecommand{\doi}[1]{doi:\discretionary{}{}{}#1}\else
  \providecommand{\doi}{doi:\discretionary{}{}{}\begingroup
  \urlstyle{rm}\Url}\fi

\bibitem[{irn({})}]{irnss}
\enquote{Indian Regional Navigation Satellite System (IRNSS),}
  \url{https://www.isro.gov.in/irnss-programme}, {}.
\newblock Accessed January 15, 2019.

\bibitem[{qzs({})}]{qzss}
\enquote{Quasi-Zenith Satellite System (QZSS),} \url{http://qzss.go.jp/en/},
  {}.
\newblock Accessed January 15, 2019.

\bibitem[{Diekelman(1998)}]{diekelman1998}
Diekelman, D., \enquote{Design guidelines for post-2000 constellations,}
  \emph{Mission Design \& Implementation of Satellite Constellations},
  Springer, 1998, pp. 11--21.

\bibitem[{Lee et~al.(2018{\natexlab{a}})Lee, Jakob, Ho, Shimizu, and
  Yoshikawa}]{lee2018}
Lee, H.~W., Jakob, P.~C., Ho, K., Shimizu, S., and Yoshikawa, S.,
  \enquote{Optimization of satellite constellation deployment strategy
  considering uncertain areas of interest,} \emph{Acta Astronautica},
  2018{\natexlab{a}}.

\bibitem[{Lutz et~al.(2012)Lutz, Werner, and Jahn}]{lutz2012}
Lutz, E., Werner, M., and Jahn, A., \emph{Satellite systems for personal and
  broadband communications}, Springer Science \& Business Media, 2012.
\newblock \urlprefix\url{https://doi.org/10.1007/978-3-642-59727-5}.

\bibitem[{Walker(1970)}]{walker1970}
Walker, J.~G., \enquote{Circular orbit patterns providing continuous whole
  earth coverage,} Tech. rep., Royal Aircraft Establishment Farnborough (United
  Kingdom), 1970.

\bibitem[{Walker(1977)}]{walker1977}
Walker, J.~G., \enquote{Continuous whole-earth coverage by circular-orbit
  satellite patterns,} Tech. rep., ROYAL AIRCRAFT ESTABLISHMENT FARNBOROUGH
  (UNITED KINGDOM), 1977.

\bibitem[{Walker(1984)}]{walker1984}
Walker, J.~G., \enquote{Satellite constellations,} \emph{Journal of the British
  Interplanetary Society}, Vol.~37, 1984, pp. 559--572.

\bibitem[{Luders(1961)}]{luders1961}
Luders, R.~D., \enquote{Satellite networks for continuous zonal coverage,}
  \emph{ARS Journal}, Vol.~31, No.~2, 1961, pp. 179--184.
\newblock \urlprefix\url{https://doi.org/10.2514/8.5422}.

\bibitem[{L{\"u}ders and Ginsberg(1974)}]{luders1974}
L{\"u}ders, R., and Ginsberg, L., \enquote{Continuous zonal coverage-a
  generalized analysis,} \emph{Mechanics and Control of Flight Conference},
  1974, p. 842.
\newblock \urlprefix\url{https://doi.org/10.2514/6.1974-842}.

\bibitem[{Beste(1978)}]{beste1978}
Beste, D.~C., \enquote{Design of satellite constellations for optimal
  continuous coverage,} \emph{IEEE Transactions on Aerospace and Electronic
  Systems}, , No.~3, 1978, pp. 466--473.
\newblock \doi{10.1109/TAES.1978.308608}.

\bibitem[{Rider(1986)}]{rider1986}
Rider, L., \enquote{Analytic design of satellite constellations for zonal earth
  coverage using inclined circular orbits,} \emph{Journal of the Astronautical
  Sciences}, Vol.~34, 1986, pp. 31--64.

\bibitem[{Ballard(1980)}]{ballard1980}
Ballard, A.~H., \enquote{Rosette constellations of earth satellites,}
  \emph{IEEE Transactions on Aerospace and Electronic Systems}, , No.~5, 1980,
  pp. 656--673.
\newblock \doi{10.1109/TAES.1980.308932}.

\bibitem[{Draim(1987)}]{draim1987}
Draim, J.~E., \enquote{A common-period four-satellite continuous global
  coverage constellation,} \emph{Journal of Guidance, Control, and Dynamics},
  Vol.~10, No.~5, 1987, pp. 492--499.
\newblock \urlprefix\url{https://doi.org/10.2514/3.20244}.

\bibitem[{Wertz(2001)}]{wertz2001}
Wertz, J.~R., \emph{Mission Geometry; Orbit and Constellation Design and
  Management: Spacecraft Orbit and Attitude Systems}, Space technology library,
  Microcosm Press, 2001.

\bibitem[{Hanson et~al.(1992)Hanson, Evans, and Turner}]{hanson1992designing}
Hanson, J.~M., Evans, M.~J., and Turner, R.~E., \enquote{Designing good partial
  coverage satellite constellations,} \emph{Journal of the Astronautical
  Sciences}, Vol.~40, No.~2, 1992, pp. 215--239.
\newblock \urlprefix\url{https://doi.org/10.2514/6.1990-2901}.

\bibitem[{Ma and Hsu(1997)}]{ma1997}
Ma, D.-M., and Hsu, W.-C., \enquote{Exact design of partial coverage satellite
  constellations over oblate earth,} \emph{Journal of Spacecraft and Rockets},
  Vol.~34, No.~1, 1997, pp. 29--35.
\newblock \urlprefix\url{https://doi.org/10.2514/6.1994-3721}.

\bibitem[{Pontani and Teofilatto(2007)}]{pontani2007}
Pontani, M., and Teofilatto, P., \enquote{Satellite constellations for
  continuous and early warning observation: A correlation-based approach,}
  \emph{Journal of guidance, control, and dynamics}, Vol.~30, No.~4, 2007, pp.
  910--921.
\newblock \urlprefix\url{https://doi.org/10.2514/1.23094}.

\bibitem[{Crossley and Williams(2000)}]{crossley2000}
Crossley, W.~A., and Williams, E.~A., \enquote{Simulated annealing and genetic
  algorithm approaches for discontinuous coverage satellite constellation
  design,} \emph{Engineering Optimization+ A35}, Vol.~32, No.~3, 2000, pp.
  353--371.
\newblock \doi{10.1080/03052150008941304},
  \urlprefix\url{https://doi.org/10.1080/03052150008941304}.

\bibitem[{Ulybyshev(2008)}]{ulybyshev2008}
Ulybyshev, Y., \enquote{Satellite constellation design for complex coverage,}
  \emph{Journal of Spacecraft and Rockets}, Vol.~45, No.~4, 2008, pp. 843--849.
\newblock \urlprefix\url{https://doi.org/10.2514/1.35369}.

\bibitem[{Dutruel-Lecohier and Mora(1998)}]{dutruel-lecohier1998}
Dutruel-Lecohier, G., and Mora, M.~B., \enquote{Orion --- A Constellation
  Mission Analysis Tool,} \emph{Mission Design {\&} Implementation of Satellite
  Constellations}, edited by J.~C. van~der Ha, Springer Netherlands, Dordrecht,
  1998, pp. 373--393.

\bibitem[{Ulybyshev(2009)}]{ulybyshev2009}
Ulybyshev, Y., \enquote{Satellite constellation design for continuous coverage:
  short historical survey, current status and new solutions,} \emph{Proceedings
  of Moscow Aviation Institute}, Vol.~13, No.~34, 2009, pp. 1--25.

\bibitem[{Mortari et~al.(2004)Mortari, Wilkins, and Bruccoleri}]{mortari2004}
Mortari, D., Wilkins, M.~P., and Bruccoleri, C., \enquote{The flower
  constellations,} \emph{Journal of Astronautical Sciences}, Vol.~52, No.~1,
  2004, pp. 107--127.

\bibitem[{Mortari and Wilkins(2008)}]{mortari2008}
Mortari, D., and Wilkins, M.~P., \enquote{Flower constellation set theory. Part
  I: Compatibility and phasing,} \emph{IEEE Transactions on Aerospace and
  Electronic Systems}, Vol.~44, No.~3, 2008, pp. 953--962.
\newblock \doi{10.1109/TAES.2008.4655355}.

\bibitem[{Wilkins and Mortari(2008)}]{wilkins2008}
Wilkins, M.~P., and Mortari, D., \enquote{Flower constellation set theory part
  ii: secondary paths and equivalency,} \emph{IEEE Transactions on Aerospace
  and Electronic Systems}, Vol.~44, No.~3, 2008, pp. 964--976.
\newblock \doi{10.1109/TAES.2008.4655356}.

\bibitem[{Lee et~al.(2018{\natexlab{b}})Lee, Ho, Shimizu, and
  Yoshikawa}]{lee2018AAS}
Lee, H.~W., Ho, K., Shimizu, S., and Yoshikawa, S., \enquote{A semi-analytical
  approach to satellite constellation design for regional coverage,}
  \emph{AAS/AIAA Astrodynamics Specialist Conference}, 2018{\natexlab{b}}.

\bibitem[{Vtipil and Newman(2012)}]{vtipil2012determining}
Vtipil, S., and Newman, B., \enquote{Determining an Earth observation repeat
  ground track orbit for an optimization methodology,} \emph{Journal of
  Spacecraft and Rockets}, Vol.~49, No.~1, 2012, pp. 157--164.
\newblock \urlprefix\url{https://doi.org/10.2514/1.A32038}.

\bibitem[{Bruccoleri(2007)}]{bruccoleri2007}
Bruccoleri, C., \enquote{Flower Constellation Optimization and Implementation,}
  Ph.D. thesis, Texas AM University, 2007.

\bibitem[{Wu et~al.(2006)Wu, Wu, and Zhu}]{wu2006design}
Wu, T., Wu, S., and Zhu, L., \enquote{Design of common track satellite
  constellations for optimal regional coverage,} \emph{2006 6th International
  Conference on ITS Telecommunications}, IEEE, 2006, pp. 1252--1255.
\newblock \doi{10.1109/ITST.2006.288854}.

\bibitem[{Avenda{\~{n}}o et~al.(2013)Avenda{\~{n}}o, Davis, and
  Mortari}]{avendano2013}
Avenda{\~{n}}o, M.~E., Davis, J.~J., and Mortari, D., \enquote{The 2-D lattice
  theory of Flower Constellations,} \emph{Celestial Mechanics and Dynamical
  Astronomy}, Vol. 116, No.~4, 2013, pp. 325--337.
\newblock \doi{10.1007/s10569-013-9493-8},
  \urlprefix\url{https://doi.org/10.1007/s10569-013-9493-8}.

\bibitem[{Chylla and Eagle(1992)}]{chylla1992}
Chylla, M.~A., and Eagle, C.~D., \enquote{Efficient computation of satellite
  visibility periods,} \emph{Spaceflight Mechanics 1992}, 1992, pp. 823--834.

\bibitem[{Alfano et~al.(1992)Alfano, Negron~Jr, and Moore}]{alfano1992}
Alfano, S., Negron~Jr, D., and Moore, J.~L., \enquote{Rapid determination of
  satellite visibility periods,} Tech. rep., AIR FORCE ACADEMY COLORADO SPRINGS
  CO, 1992.

\bibitem[{Han et~al.(2017)Han, Gao, and Sun}]{han2017}
Han, C., Gao, X., and Sun, X., \enquote{Rapid satellite-to-site visibility
  determination based on self-adaptive interpolation technique,} \emph{Science
  China Technological Sciences}, Vol.~60, No.~2, 2017, pp. 264--270.
\newblock \urlprefix\url{https://doi.org/10.1007/s11431-016-0513-8}.

\bibitem[{stk({})}]{stk}
\enquote{Systems Tool Kit Help Guide,} \url{http://help.agi.com/stk/}, {}.
\newblock Accessed April 28, 2019.

\bibitem[{Gray(2006)}]{gray2006toeplitz}
Gray, R.~M., \enquote{Toeplitz and circulant matrices: A review,}
  \emph{Foundations and Trends{\textregistered} in Communications and
  Information Theory}, Vol.~2, No.~3, 2006, pp. 155--239.
\newblock \doi{10.1561/0100000006}.

\bibitem[{Oppenheim and Schafer(2009)}]{dsp}
Oppenheim, A.~V., and Schafer, R.~W., \emph{Discrete-Time Signal Processing},
  3\textsuperscript{rd} ed., Prentice Hall Press, USA, 2009.

\bibitem[{Gurobi~Optimization(2016)}]{gurobi}
Gurobi~Optimization, I., \enquote{Gurobi Optimizer Reference Manual,} , 2016.
\newblock \urlprefix\url{http://www.gurobi.com}.

\bibitem[{Lee et~al.(2016)Lee, Wu, and Mortari}]{lee2016satellite}
Lee, S., Wu, Y., and Mortari, D., \enquote{Satellite constellation design for
  telecommunication in Antarctica,} \emph{International Journal of Satellite
  Communications and Networking}, Vol.~34, No.~6, 2016, pp. 725--737.
\newblock \urlprefix\url{https://doi.org/10.1002/sat.1128}.

\bibitem[{Wessel and Smith(2014)}]{antarctica}
Wessel, P., and Smith, W. H.~F., \enquote{Shoreline Boundary between Antarctic
  Grounding Line and the Ocean, 2014 (Full-resolution),} , 2014.
\newblock \urlprefix\url{http://purl.stanford.edu/zt046tb6131}.

\bibitem[{Chambers et~al.(2007)Chambers, Asner, Morton, Anderson, Saatchi,
  Esp{\'\i}rito-Santo, Palace, and Souza~Jr}]{chambers2007regional}
Chambers, J.~Q., Asner, G.~P., Morton, D.~C., Anderson, L.~O., Saatchi, S.~S.,
  Esp{\'\i}rito-Santo, F.~D., Palace, M., and Souza~Jr, C., \enquote{Regional
  ecosystem structure and function: ecological insights from remote sensing of
  tropical forests,} \emph{Trends in Ecology \& Evolution}, Vol.~22, No.~8,
  2007, pp. 414--423.
\newblock \urlprefix\url{https://doi.org/10.1016/j.tree.2007.05.001}.

\bibitem[{Rientjes et~al.(2013)Rientjes, Haile, and
  Fenta}]{rientjes2013diurnal}
Rientjes, T., Haile, A.~T., and Fenta, A.~A., \enquote{Diurnal rainfall
  variability over the Upper Blue Nile Basin: A remote sensing based approach,}
  \emph{International journal of applied earth observation and geoinformation},
  Vol.~21, 2013, pp. 311--325.
\newblock \urlprefix\url{https://doi.org/10.1016/j.jag.2012.07.009}.

\bibitem[{the({})}]{theworldbasin}
\enquote{Major River Basins Of The World,}
  \url{https://datacatalog.worldbank.org/dataset/major-river-basins-world}, {}.
\newblock Accessed January 11, 2020.

\bibitem[{Kim(1997)}]{kim1997theory}
Kim, M., \enquote{Theory of satellite ground-track crossovers,} \emph{Journal
  of Geodesy}, Vol.~71, No.~12, 1997, pp. 749--767.
\newblock \urlprefix\url{https://doi.org/10.1007/s001900050141}.

\end{thebibliography}

\newpage

\end{document}